\documentclass[11pt,twoside]{article}
\usepackage{latexsym}
\usepackage{envmath}
\usepackage{amssymb,amsbsy,amsmath,amsfonts,amssymb,amscd}
\usepackage{graphicx}
\usepackage{xcolor}
\usepackage{hyperref}
\usepackage{enumitem}
\usepackage[
backend=biber,
style=numeric,
sorting=nty,maxbibnames=99,maxcitenames=99
]{biblatex}
\renewbibmacro*{in:}{%
  \ifboolexpr{
    test {\ifentrytype{article}}
  }
    {} 
    {\printtext{\bibstring{in}\intitlepunct}}
}
\numberwithin{equation}{section}
\newcommand{\commentout}[1]{}
\newcommand{\R}{\mathbb{R}}
\newcommand{\N}{\mathbb{N}}

\newcommand {\e}  {\varepsilon}

\newcommand {\Chi} {{\bf \raise 2pt \hbox{$\chi$}} }
\newcommand{\ued}{u_{\varepsilon,\delta}}
\newcommand {\f}   {\frac}
\newcommand {\p}   {\partial}
\newcommand{\fer}{\eqref}
\newcommand{\dis}{\displaystyle}
\newcommand {\proof} {\noindent {\bf Proof}. }
\newcommand{\beq}{\begin{equation}}
\newcommand{\eeq}{\end{equation}}
\newcommand{\bea} {\begin{array}{rl}}
\newcommand{\eea} {\end{array}}
\newcommand{\bepa}{\left\{ \begin{array}{l}}
\newcommand{\eepa} {\end{array}\right.}
\newtheorem{theorem}{Theorem}[section]
\newtheorem{lemma}[theorem]{Lemma}

\newtheorem{remark}[theorem]{Remark}
\newtheorem{proposition}[theorem]{Proposition}
\newtheorem{corollary}[theorem]{Corollary}
\newcommand{\qed}{{ \hfill
                       {\unskip\kern 6pt\penalty 500 \raise -2pt\hbox{\vrule\vbox to 6pt{\hrule width 6pt
                       \vfill\hrule}\vrule} \par}   }}
\title{Concentration in selection-mutation models: error estimates and asymptotic expansions}

\author{
Caroline Guinet \thanks{Institut de Math\'ematiques de Toulouse; UMR 5219, Universit\'e de Toulouse; CNRS, UPS IMT, F-31062 Toulouse Cedex 9, France; E-mail: caroline.guinet@math.univ-toulouse.fr} \and
Sepideh Mirrahimi \thanks{Institut de Math\'ematiques de Toulouse; UMR 5219, Universit\'e de Toulouse; CNRS, UPS IMT, F-31062 Toulouse Cedex 9, France; E-mail: sepideh.mirrahimi@math.univ-toulouse.fr} 
\and
Jean-Michel Roquejoffre  \footnotemark[2]\thanks{Institut de Math\'ematiques de Toulouse; UMR 5219, Universit\'e de Toulouse; CNRS, UPS IMT, F-31062 Toulouse Cedex 9, France; E-mail: jean-michel.roquejoffre@math.univ-toulouse.fr}}

\date{\today}

\begin{document}
\maketitle
\pagestyle{plain}

\pagenumbering{arabic}

\begin{abstract}
In this paper, we study an integro-differential equation which describes the evolutionary dynamics of a population structured by a phenotypic trait. This population undergoes asexual reproduction, competition, selection, and mutation. We provide an asymptotic analysis of the model, assuming that the mutations have small effects.  A standard approach for the analysis of the qualitative properties of the solutions of such an equation is to apply a logarithmic transformation, which yields a Hamilton–Jacobi equation with constraint. When the reproduction term is a concave function of the trait, it has been established that the solution is classical. We rigorously derive a first-order asymptotic expansion of the solution. This expansion allows us to approximate the moments of the phenotypic density. This result establishes a connection between the approximations of the phenotypic density obtained via the Hamilton-Jacobi approach and relevant biological quantities, which are more suitable from a modeling perspective.
\end{abstract}
\textbf{AMS Class. No:} 35B25, 35C20, 35F21, 35Q92, 92D15.
\bigskip
\tableofcontents


\section{Introduction}

\subsection{Model and main question}
In this paper, we study the asymptotic behavior of the following parabolic Lotka-Volterra type integro-differential equation 
\begin{equation} \label{para}
\e\p_t n_\e - \e^2 \Delta n_\e=  n_\e R \big(x, I_\e(t) \big), \ t>0, \; x \in \R^d,
\end{equation}
with a nonlinearity driven by the integral term
\begin{equation} \label{paraI}
I_\e(t) = \int_{\R^d}  \psi(x) n_\e(t,x) \mathrm{d}x,
\end{equation}
with $\psi$ a given smooth positive function and $n_{\e,0}(x)$ the initial condition. Our objective is to estimate $n_\e$ in the limit of small $\e$. Under some monotonicity and concavity assumptions, it has been shown in \cite{GB.BP:08,GB.SM.BP:09,AL.SM.BP:10,SM.JR:16} that one can write $$ n_\e = \frac{r}{\e^{\frac{d}{2}}}\exp \Big( \frac{u+o(1)}{\e}\Big).$$
The objective here is to extend the approximation to the next order and prove that $$ n_\e = \frac{r}{\e^{\frac{d}{2}}}\exp \Big( \frac{u+ \e v+ o(\e)}{\e}\Big). $$
\subsection{Motivation and the state of the art}
The equation \eqref{para} is known as the selection-mutation model. It describes the evolutionary dynamics of a population structured by a phenotypic trait. The population undergoes asexual reproduction, competition, selection and mutations. The solution $n_{\e}(t,x)$ stands for the density of individuals with trait $x \in  \R^d$, at time $t.$ The mutations are represented  by a Laplace term with a small coefficient, assuming in this way that the mutations have small effects. The individuals with trait $x$ at time $t$ reproduces with rate $R(x,I_{\e}(t))$. This rate depends on the environmental feedback $I_\e(t).$ This integral term can represent the total consumption of a nutrient. The dependency in $I_\e(t)$ of the growth rate $R(x,I_{\e}(t))$ models the competition in the population. Indeed, we will choose $R$ to be strictly decreasing with respect to $I.$ A variant of this model was first introduced by Kimura \cite{MK:65}.  It was later derived from a stochastic individual-based model with small mutational effects but a high mutation rate in \cite{NC.RF.SM:06} and \cite{NC.RF.SM:08}. \\\\
This class of models has been widely investigated in the literature using diverse approaches. See, for instance, \cite{OD:04} for an approach based on the stability analysis of differential systems or \cite{NC.RF.SM:06,NC.RF.SM:08,champagnat:tel-00091929} for stochastic approaches considering individual-based models, and \cite{GB.BP:08,GB.SM.BP:09,AL.SM.BP:10,LD.PJ.SM.GR:08,PJ.GR:09,Raoulphd,FH.FL.LR} for the study of integro-differential models.
\\\\
 Here, we focus on a method involving Hamilton-Jacobi equations, considering small mutational effects. This method has been first developed in \cite{OD.PJ.SM.BP:05,GB.BP:08,GB.SM.BP:09,AL.SM.BP:10}. The aim is to study the qualitative properties of the density $n_\e$ when $\e$ approaches $0$. The solution to \eqref{para} typically concentrates around a dominant phenotypical trait $\overline{x}(t)$ that evolves in time, i.e. $$n_\e  \rightharpoonup_{\e \rightarrow 0} \rho(t)\delta(x-{\overline{x}(t)}),$$
where $\delta(x-{\overline{x}(t)})$ is a Dirac mass centered in $\overline x(t).$ 
Additionally, a closely related approach was used for reaction–diffusion equations in the context of the ‘geometric optics’ approximation: Freidlin \cite{MF:85,MFb:85} employed a probabilistic method, while Barles, Evans, and Souganidis \cite{LE.PS:89,GB.LE.PS:90} used viscosity solutions. The first step in this approach is to perform the Hopf-Cole transformation: 
\begin{equation}
    \label{HC}
    n_\e(t,x) = \exp \Big( \frac{u_\e(t,x)}{\e} \Big).
\end{equation}
The function $u_\e$ satisfies the following equation 
\begin{equation}
    \label{eq-u_e}
    \partial_t u_\e -\e \Delta u_\e = |\nabla u_\e|^2 + R(x,I_\e(t)).
\end{equation}
Relying on an analysis of \eqref{eq-u_e}, one can prove \cite{GB.BP:08,GB.SM.BP:09,AL.SM.BP:10} the convergence along subsequences of $u_\e$ to $u$ a viscosity solution of the following non classical Hamilton-Jacobi equation with constraint 
\begin{System}
\label{HJ3}
    \partial_t u =|\nabla u |^2 + R(x,I(t)),\ t >0, \ x \in \R^d,\\
    \max_{x \in \R^d} u(t,x) = 0, \\
    u(0,x) = u^0.
\end{System}
It can be noticed that the function $I(t)$ is now an unknown of problem 
\eqref{HJ3}, it will be determined by the constraint on $u.$ We need uniqueness of equation $(u,I)$ of the latter solution to prove strong convergence of $u_\e$ to $u$. This uniqueness problem has been studied in different frameworks before. It was treated in \cite{GB.BP:08} considering a particular form of $R(x,I)$. Later, it was proved in \cite{SM.JR:16} in the concave framework, i.e. assuming $R$ and $u_0$ to be strictly concave with respect to $x.$ It was also proved later, with relaxed assumptions, in \cite{VC.KL:19}, within the class of BV functions. The uniqueness property implies that the convergence along subsequences is indeed a strong convergence of the whole sequence.\\\\
The goal of this paper is to derive an asymptotic expansion for the solution $u_\e$ and for the integral term $I_\e$ when $\e$ approaches $0.$ The interest of this work lies (i) in the mathematical exploration of a singular limit, and (ii) in providing a rigorous access to quantities of interest in theoretical biology. Formal asymptotic expansions were already used by Figueroa Iglesias, Mirrahimi in \cite{SF.SM:18} and by Gandon, Mirrahimi in \cite{SG.SM:18}, respectively for time-periodic and space heterogeneous environments. This enables them to derive asymptotic expansions for the moments of the phenotypic distribution, which are more directly measurable in biology. The results from \cite{SF.SM:18} and \cite{SG.SM:18} demonstrate that these approximations yield more informative conclusions for applications.
It is natural to focus on this simple model, as its study may serve as a basis for more elaborated models incorporating additional effects.
\subsection{Assumptions}
The assumptions  we are stating below are the same as in \cite{AL.SM.BP:10}, where the authors noticed that this set of assumptions allowed them to work with smooth solutions, thus going quite far in the study of \fer{para}. We believe that the results that we will prove under these assumptions would certainly be false if some of those assumptions  were removed. 
\begin{itemize}
\item {\bf Assumptions on $\psi$.}  The function $\psi$ is chosen so that
\begin{equation} \label{aspsi}
0<\psi_m\leq \psi \leq \psi_M < \infty, \qquad \psi \in W^{3,\infty}(\R^d)
.
\end{equation}
\item{\bf Assumptions on $R$.} We choose $R\in C^2$, and we suppose  that  there is $ I_M>0$  such that (fixing the origin in $x$ appropriately)
\begin{equation} \label{asrmax}
\max_{x \in \R^d} R(x,I_M) = 0 = R(0,I_M).
\end{equation}
This assumption implies that the resources are limited and that the growth rate $R$ is nonpositive for every individual in a population of weighted size $I_M.$	We also assume that
\begin{equation} \label{asr}
-\underline{K}_1 |x|^2 \leq R(x,I) \leq \overline{K}_0 -\overline{K}_1 |x|^2,  \qquad \text{for }\;  0\leq I \leq I_M,
\end{equation}
\begin{align}
\label{asrD2}
&- 2\underline{K}_1 \leq D^2 R(x,I) \leq - 2\overline{K}_1 < 0    \text{ as symmetric matrices,  } \nonumber \\& D^k R(\cdot,I) \in L^\infty(\R^d),\ \text{for $k \in \{3,...,6\}$ and for } 0\leq I\leq I_M.
\end{align} We also assume that
\begin{equation} \label{asrDi}
- \underline{K}_2\leq \dis{\frac{\p R}{\p I} \leq - \overline{K}_2}.
\end{equation}
These assumptions imply in particular that for all $I$, $R$ has a single maximum point, that is there is only one optimal trait for all values of $I$. However, the concavity of 
$R$ is a stronger technical assumption which guarantees, together with a concavity assumption on the initial data, that the solution to \eqref{HJ3} remains smooth and concave for all times. The monotonicity assumption models the competition between individuals, it means that the growth rate decreases as the size of the weighted population increases.\\\\
We also assume that
\beq
\label{asr23}
| \f{\p^2 R}{\p I \p x_i}(x,I)| + | \f{\p^3 R}{\p I \p x_{i} \p x_j}(x,I)| + | \f{\p^4 R}{\p I \p x_{i} \p x_j \p x_k}(x,I)| \leq K_3,  \ 
\eeq
$\text{for  $0\leq I \leq I_M$, and $i,\,j,\ k =1,2,\cdots,d$}.$
\item{\bf Assumptions on $n(0,.)$.}
To study the qualitative behavior of $n_\e$ we use, as is by now classical,
 the Hopf-Cole transformation
\beq
\label{Hopf}
n_\e =\exp \left( \f {u_\e} {\e} \right).
\eeq
The initial datum $n(0,.):=n_\e^0$ will be chosen to satisfy some compatibility conditions  with the assumptions on $R$ and $\psi$. Namely, we require that there is  $I^0$ such that
\begin{equation} \label{asI}
0<I^0\leq I_\e(0) := \int_{\R^d} \psi(x) n_\e^0(x) dx< I_M,
\end{equation}
and we can write
$$
n_\e^0= e^{\f{u_\e^0}{\e}}=\f{r}{ \e^{\f d 2}}\, e^{\f{u^0}{\e}},  \ \text{with }\;  u^0 \in C^2(\R^d), \ \text{and  } \max_{x\in \R^d} u^0(x)=0,
$$
which implies that
$$ 
u_\e^0 = u^0  + \e \log \left( \f{r}{\e^{\f d 2}} \right).
$$
As for $u^0$ we assume the existence of positive constants $\underline{L}_0, \overline{L}_0,\underline{L}_1, \overline{L}_1$ such that
\begin{equation} \label{asu}
-\underline{L}_0 -\underline{L}_1 |x|^2 \leq u^0(x)\leq  \overline{L}_0 - \overline{L}_1 |x|^2 ,
\end{equation}
\begin{equation} \label{asuD2}
-2\underline{L}_1 \leq D^2u^0  \leq - 2\overline{L}_1 .
\end{equation}
We also assume that 
\beq \label{extinct} \frac{1}{\e}\Big(  \int_{\R^d} \psi(x) R(x,I_\e(0))n_\e^0(x) dx\Big)_{-} = o(1), \ \text{as } \e \text{ goes to } 0.\eeq In terms of modeling, this assumption guarantees the survival of the population.\\\\
We also need the following technical assumptions 
\begin{equation} \label{asuD3}
D^k u^0 \in L^\infty(\R^d),\ \text{for $k\in \{3,...6\},$}
\end{equation}
\begin{align} \label{asuIni}
&n_\e^0(x) \underset{\e\to 0}{\longrightarrow}  \bar{\rho}^0 \; \delta \big(x-x_0\big) \text{ weakly in the sense of measures} \nonumber,\\ &\text{with $R(x_0,I_0)=0$ and $I_0=  \psi(x_0)\bar{\rho}^0$.}
\end{align}
\end{itemize}
\subsection{Notations}

\begin{itemize}
\item To describe the asymptotic behavior of a function $f_\e$, when $\e$ approaches $0$, we will use the following notations. We will say that $f_\e = O(g_\e)$ if and only if there exists $C>0$ independent of $\e$, such that for $\e$ sufficiently small, we have $|f_\e| \leq C |g_\e|,$ uniformly in time and space. We also say that $f_\e = o(g_\e)$ if and only if there exists a function $\alpha_\e$, such that for $\e$ small enough we have
$f_\e=\alpha_\e g_\e$, and $\alpha_\e$ tends uniformly to $0$ when $\e$ goes to $0$.
\item In this paper, all constants 
$C$ are independent of $\e$, may eventually depend on the time interval length $T$, and can change from one line to another.
\item For all $a,b \in \R,$ we define the interval $e(a,b)$ as follows 
$$ e(a,b) := (\min(a,b), \max(a,b)).$$
\item We introduce notations corresponding to the moments of the phenotypic distribution. Here we choose 
 $\psi \equiv 1$ and we define $q_\e:= \frac{n_\e}{
I_\e}$ to be the phenotypical distribution. We then introduce notations for the first moment and the central moments of order $k$: 
\begin{align*}
    M_{1,\e}(t) &:= \int_{\R^d} x q_\e(t,x) \mathrm{d}x, \\
    M_{k,\e}^c(t) &:= \int_{\R^d} (x-M_{1,\e}(t))^k q_\e(t,x) \mathrm{d}x.
\end{align*}
We will use these notations to examine the approximations of moments.
\end{itemize}
\subsection{Previous results}
Under the above assumptions the following theorem is proved in \cite{AL.SM.BP:10}.

\begin{theorem}\cite{AL.SM.BP:10}
\label{th:HJ}
Under assumptions \fer{aspsi}--\fer{asuIni}, as $\e \to 0$ and along subsequences, $(u_\e)_\e$ converges locally uniformly to
 $u\in L^\infty_{\mathrm loc} \big(R^+; W^{3,\infty}_{\mathrm loc}(\R^d) \big)\cap W^{1,\infty}_{\mathrm loc} \big(R^+; L^{\infty}_{\mathrm loc}(\R^d) \big)$ and $(I_\e)_\e$ converges almost every where to $I  \in W^{1,\infty}(\R^+)$, where $u$ is a viscosity solution of the Hamilton-Jacobi equation with constraint \eqref{HJ3}.
Moreover, for all $t\geq 0$, the maximum of $u_\e(t,\cdot)$ is attained at a unique point $ x_\e(t)\in \R^d$, $ x_\e (\cdot) \in W^{1,\infty}([0,T]; \R^d)$ uniformly in $\e\leq \e_0$, for all $T>0$ and  some $\e_0>0$. As $\e\to 0$, $(x_\e(t))_\e$ converges locally uniformly along subsequences to $ \overline x(t) \in W^{1,\infty}(\R^+; \R^d)$, the unique maximum point of $u(t,\cdot)$, which satisfies
\begin{align}\label{eq.cano}
\dot{ \overline x }(t)= \left( -D^2u \big(t, \overline x(t)\big) \right)^{-1} \nabla_x R\big(\overline x(t),\bar I(t) \big), \quad \overline x(0)={x}_0,
\end{align}
and 
\beq
\label{R=0}
R \left( \overline x(t) , I(t) \right) =0,\qquad I(t)\leq I_M.
\eeq
\end{theorem}
\begin{corollary}
As $\e\to 0$ and along subsequences, we have
$$
n_\e(t,x) \underset{\e\to 0}{\longrightarrow}  \bar \rho(t) \; \delta \big(x-\overline x (t)\big),\qquad \text{weakly in the sense of measures},
$$
with 
$$
\rho(t)= \f{I(t) }{\psi(\overline x(t))}.
$$
\end{corollary}
A natural question is therefore to understand whether the convergence in Theorem \ref{th:HJ} can be improved, i.e. (i) whether it holds in a stronger sense than 
along subsequences, (ii) whether an asymptotic expansion of the solution can be provided. The following uniqueness theorem, proved in \cite{SM.JR:16}, answers the question (i).
\begin{theorem} \cite{SM.JR:16}
\label{thm:uniq}
Under assumptions \fer{aspsi}--\fer{asuIni}, the Hamilton-Jacobi equation with constraint \eqref{HJ3} has a unique solution $(u,I) \in L^\infty_{\mathrm loc} \big(\R^+; W^{3,\infty}_{\mathrm loc}(\R^d) \big)\cap C^{1}\big(\R^+\times \R^d \big) \times C^{1}(\R^+)$ and $\nabla u \in C^{1}\big(\R^+\times \R^d \big).$
\end{theorem}
The uniqueness result was established  in the general framework and without the concavity assumption by Calvez and Lam in \cite{VC.KL:19}. The uniqueness replaces convergence along subsequences by convergence of the full family, and allows for error estimates.  In this paper, we will provide an asymptotic expansion of the solution (ii).
\subsection{Main results} 
The next Theorem proves first order asymptotic expansions for the solution $u_\e$ of \eqref{eq-u_e}, its maximum $x_\e$ and the quantity $I_\e$.
\begin{theorem}
\label{thm:approx}
We assume \fer{aspsi}--\fer{asuIni}.
Let $n_\e$ be the solution of \fer{para} and $u_\e$ be defined by \fer{Hopf}. There exist continuous functions  $J:\R^+\to \R$, $v \in L^\infty_t W^{2,\infty}_x([0,T]\times \R^d)$, and a constant $C(T)=C>0$ independant of $\e$, such that we have the following expansions:
\beq
\label{expI}
    \|I_\e - I - \e J \|_{L^\infty([0,T])}\leq C \e^2,
\eeq 
\beq 
\label{expx}
    \|x_\e - \overline x - \e y  \|_{L^\infty([0,T])} \leq C\e^2,\eeq
\beq
\label{expu}
   \| u_\e - u - \e v  -\e\log \Big(\frac{r}{\e^{d/2}}\Big) \|_{L^\infty_t W^{2,\infty}_x ([0,T]\times \R^d)} \leq C \e^2,
\eeq
where $x_\e(t)$ is the maximum point of $u_\e$:
 $$
 \max_{x\in \R^d}u_\e(t,x) =u_\e(t,x_\e(t)).
 $$
\end{theorem}
This proof relies
on the methods of Mirrahimi and Roquejoffre in \cite{SM.JR:16}, in addition to the results of \cite{AL.SM.BP:10} within the concave framework. They noticed that the solution $u$ is smooth, strictly concave and that the unique maximum of $u$ evolves according to an ODE. Moreover, we exploit the equivalence proved in \cite{SM.JR:16} between the constrained Hamilton-Jacobi problem \eqref{HJ3} and the ODE-PDE formulation \eqref{system-e}. The idea is then to first control $I_\e$ using the Laplace's integration method and then compare $u_\e$ to $u$ using the PDE-ODE formulation. To this end, we will use regularity estimates, maximum principle, and the method of characteristics. Using the asymptotic expansions from Theorem \ref{thm:approx}, we obtain first order asymptotic expansions for the moments of the phenotypic distribution.
\begin{theorem}
\label{thm:moments}
    We assume \fer{aspsi}--\fer{asuIni}, and  $\psi \equiv 1$. For all $k\geq 1,$ there exists a constant $C_k(T) = C_k>0,$ independant of $\e$, such that we have
    \begin{align}
       \| M_{1,\e}-\overline{x} -\e M_1\|_{L^\infty([0,T])}&\leq      C_1 \e \sqrt{\e}, \label{m1}\\
        \| M_{2,\e}^c -  \e M_2\|_{L^\infty([0,T])} & \leq C_2 \e^2 ,\label{m2}\\
        \| M_{2k,\e}^c - M_{2k}\e^k\|_{L^\infty([0,T])} & \leq C_{2k}\e^k \sqrt{\e}, \label{m2k} \\ 
        \| M_{2k+1,\e}^c- M_{2k+1}\e^{k+1}\|_{L^\infty([0,T])} & \leq C_{2k+1}\e^{k+1}, \label{m2k1}
    \end{align}
    with 
\begin{align} \label{M11} M_1(t) = & \ \frac{1}{2}\big(-D^2 u(t,\overline{x}(t))\big)^{-{2}} D^3 u(t,\overline{x}(t)) +  \big(- D^2 u(t,\overline{x}(t))\big)^{-1}\nabla v (t,\overline{x}(t)),\\ M_2(t) =& \ \big(-D^2u(t,\overline{x}(t))\big)^{-1}, \nonumber \\M_{2k}(t) =& \ \frac{(2k)!}{2^k k!}M_2(t)^k, \nonumber \\
M_{2k+1}(t) =& \  \frac{k}{3}\cdot\frac{\big(2(k+1)\big)!}{2^{k+1}(k+1)!}D^3u(t,\overline{x}(t))M_2(t)^{k+2} \nonumber. \end{align} \end{theorem}
Theorem \ref{thm:moments} follows directly from the asymptotic expansion of Theorem  \ref{thm:approx}, combined with the Laplace’s method of integration. The proof of Theorem \ref{thm:moments} is provided in the appendix \ref{A_mom}.
\begin{remark}
We notice that the mean phenotypic trait is close to the dominant phenotypic trait $\overline x(t)$ and that the phenotypical variance is of order $\e (-D^2 u(t,\overline x(t)))^{-1}.$ Using these properties, the canonical equation \eqref{eq.cano} and after a change of variable $t \rightarrow \e t$, we obtain that the derivative of the average trait is equal to the phenotypical variance times the selection gradient. This is a generalization of Lande's equation in quantitative genetics \cite{RL:79, RL.SA:83}. This equation has been derived assuming that the phenotypic distribution is a normal distribution with constant variance. Here we relax such an assumption and capture the shape of the phenotypic distribution which is not necessarily of Gaussian type. We also track the dynamics of the phenotypic variance which may vary significantly in long times.
\end{remark}
The paper is organized as follows. In Section 2 we gather some useful facts about $u_{\varepsilon}$. All the results of this section were stated in \cite{AL.SM.BP:10}, but, due to their importance for the sequel, we reprove the concavity properties of $u_\e$. The asymptotic expansion process starts in Section 3, where we find an approximation of $I_\varepsilon(t)$. In Section 4 we prove that $u_\e, \ I_\e,$ and $u, \ I$ differ by an order at most $\e$. This will allow for a complete expansion of $I_\varepsilon$ and $u_\e$. Theorem \ref{thm:approx} is proved in Section 5.
\section{\texorpdfstring{Known features of $u_\e$ and $I_\e$}{}}
\label{sec:pre}
The first quantity to estimate is $I_\e$. Once this is under control, one may derive estimates on $u_\e$, the main one being its strict concavity. Notice here that some constants are meaningful, and that one should carefully keep track of their respective values. 
\begin{theorem}
\label{t2.1}
There exists $I_m>0$ such that
\beq
\label{bIe}
0< I_m \leq I_\e (t) \leq I_M+C\e^2 \quad \text{a.e.}
\eeq
\end{theorem}
This is a nontrivial property; see \cite{AL.SM.BP:10},  Section 4.
\begin{theorem}
\label{t2.2}
We have the following estimates on $u_\e$ for all $t\geq 0$,
\beq
  \label{eq.le.u}
-\underline{L}_0 -\underline{M}_1 |x|^2  -2\e d\underline{M}_1t \leq  u_\e(t,x) \leq \overline{L}_0 -\overline{M}_1 |x|^2 + \overline{K}_0 t,
\eeq
\beq 
\label{eq.D2u}
 -2 \underline M_1\leq D^2u_\e(t,x)  \leq -2\overline{M}_1,
\eeq
\beq
\label{nD3u}
\|D^k u_\e(t,\cdot) \|_{L^\infty(\R^d)} \leq C(T), \ \text{for $t\in [0,T]$}, \ k = \{3,...,6\}
\eeq
with $\underline{M}_1 := \max \big(\underline{L}_1, \frac{\sqrt{ \underline{K}_1}}{2} \big),$ and $\overline{M}_1 :=\min \big(\overline{L}_1, \frac{\sqrt{ \overline{K}_1}}{2} \big).$
Equation \eqref{eq.D2u} should, once again, be understood in the sense of symmetric matrices.
\end{theorem}
\begin{remark}
The estimates \eqref{eq.le.u} and \eqref{eq.D2u} were already proven by Lorz, Mirrahimi and Perthame in \cite{AL.SM.BP:10}. For completeness, we reprove these results in detail and establish estimates for the higher-order derivatives of $u_\e$. The estimate \eqref{eq.D2u} implies that $u_\e$ is strictly concave. Therefore, $u_\e(t,\cdot)$ has a unique maximum that we denote $x_\e(t)$ for the rest of the paper, i.e.
\beq \label{xe}\textrm{argmax} \ u_\e(t,\cdot) = x_\e(t).\eeq     
\end{remark} 
Before we demonstrate Theorem \ref{t2.2}, we prove the following Lemma \ref{lemme0}: 
\begin{lemma}
\label{lemme0}
    Let fix $(t_0,x_0) \in \R^+ \times \R^d$.
    We define $\overline w$ and $\underline{w}$, as $ \overline w(t_0,x_0) := \max_{\vert \xi \vert =1} \partial_{\xi\xi}^2 u(t_0,x_0)$ and $ \underline w(t_0,x_0) := \min_{\vert \xi \vert =1} \partial_{\xi\xi}^2 u(t_0,x_0)$. \\\\ Then, there exist $ \overline \xi_0, \underline \xi_0 \in \R^d $ such that $|\overline \xi_0 | = |\underline \xi_0 | = 1$, $$\overline w(t_0,x_0) = \partial_{\overline \xi_0 \overline \xi_0}^2 u(t_0,x_0), \quad D^2u(t_0,x_0). \overline \xi_0 = \overline w(t_0,x_0) \overline \xi_0,$$ and
    $$\underline w(t_0,x_0) = \partial^2_{\underline \xi_0 \underline \xi_0}u(t_0,x_0), \quad D^2u(t_0,x_0). \underline \xi_0 = \underline w(t_0,x_0) \underline \xi_0. $$
\end{lemma}
\textbf{Proof.}
We prove the result for $\overline w$, the one for $\underline w$ can be proven similarly. We define $H:= D^2u(t_0,x_0)$. \\
The matrix $H$ is real symmetric, according to the spectral theorem there exists $P \in O_d(\R)$, an orthogonal matrix, and $D = \text{diag}(\lambda_1,...,\lambda_d)$ such that $H = PDP^T$.\\
Thus, we have: $$ \forall i \in \{1,...,d\}, \quad \lambda_i = \sum_{k=1}^d \sum_{j=1}^d P_{ki}P_{ji}\partial_{kj}^2u(t_0,x_0).$$
Moreover, if $\partial_x = \sum_{k=1}^d x_k \partial_k$, then we have:
$$ \partial_{xx}^2 u(t_0,x_0) = \sum_{k=1}^d \sum_{j=1}^d x_kx_j \partial_{kj}^2 u(t_0,x_0).
$$
If $P_i$ denotes the i-th column of $P$ then for all $i \in \{1,...,d\}$, $\lambda_i$ satisfies:  $$\lambda_i = \partial_{P_iP_i}^2 u(t_0,x_0). $$
Let us show there exists $ i_0 \in \{1,...,d\}$, such that $\lambda_{i_0} = \partial_{\overline \xi_0 \overline \xi_0}^2 u(t_0,x_0)$.
\\By contradiction, we assume for all $ i \in \{1,...,d\}, \partial_{P_iP_i}^2u(t_0,x_0) \neq \partial_{\overline \xi_0 \overline \xi_0}^2 u(t_0,x_0).$ \\\\
We know by definition that $ \max_{\vert \xi \vert =1} \partial_{\xi\xi}^2 u = \partial_{\overline\xi_0\overline\xi_0}^2 u(t_0,x_0) $ thus, for all $ i \in \{1,...,d\}, \partial_{P_iP_i}^2u(t_0,x_0) < \partial_{\overline\xi_0\overline\xi_0}^2 u(t_0,x_0).$ We compute:
\begin{align*}
\partial_{\overline\xi_0\overline\xi_0}^2 u(t_0,x_0) &= (P\overline\xi_0)^T D (P\xi_0)\\
    &= \sum_{i=1}^d \lambda_i (P\overline\xi_0)^2_i < \left ( \sum_{i=1}^d (P \overline \xi_0)^2\right ) \partial_{\overline\xi_0\overline\xi_0}^2 u(t_0,x_0).
\end{align*}
We use that $P$ is an orthogonal matrix, ie. $\sum_{i=1}^d (P \overline\xi_0)^2 =1 $.
This is a contradiction, thus there exists $ i_0 \in \{1,...,d\}$ such as $\lambda_{i_0} = \partial_{\overline\xi_0 \overline\xi_0}^2 u(t_0,x_0) $. Moreover $ P_{i_0}$ is an eigenvector of $H$, of eigenvalue $\partial_{\overline\xi_0\overline\xi_0}^2 u(t_0,x_0)$. Therefore, we conclude that: 
  $$D^2u(t_0,x_0).\overline\xi_0=\nabla u_{\overline\xi_0} = \partial_{\overline\xi_0\overline\xi_0}^2 u(t_0,x_0)\overline\xi_0= \overline w(t_0,x_0)\overline\xi_0.$$
\begin{flushright}
    $\square$
\end{flushright}
\textbf{Proof of Theorem \ref{t2.2}}.\\\\
\textbf{(i) Let us prove the first estimate \fer{eq.le.u}.} \\\\
We define $a(t,x) := -\underline{L}_0 -\underline{M}_1 |x|^2  -2 d \e \underline{M}_1t$  and  $b(t,x):= \overline{L}_0 -\overline{M}_1 |x|^2 + \overline{K}_0 t$.
\begin{itemize}[leftmargin=*]
    \item Let us show that $a \leq u_\e$ :
\\We use the definition of $a$, and the assumption \fer{asr}  on $R$, to obtain: 
\begin{align*}
    \partial_t a -\e\Delta a - \vert \nabla a\vert ^2 -R(x,I_\e(t)) &= -2d\e\underline{M}_1 + 2d\e\underline M_1 - 4\underline M_1^2\vert x\vert^2 -R(x,I_\e(t)) \\
    & \leq - 4\underline M_1^2 \vert x\vert^2 + \underline K_1 \vert x\vert^2 \leq 0.
\end{align*}
with $a(0,x) =  -\underline{L}_0 -\underline{M}_1 |x|^2 \leq  u_\e(0,x),$ for $ x \in \R^d.$
Thus, $a$ is a sub-solution of \fer{eq-u_e}, we conclude using the maximum principle.
    \item Let us show that $u_\e \leq b$ :
    \\
We use the definition of $v$, and the assumption \fer{asr}  on $R$, to obtain: 
\begin{align*}
\partial_t b -\e\Delta b - \vert \nabla b\vert ^2 -R(x,I_\e(t))&= \overline K_0  +2d\e\overline{M}_1 -4 \overline M_1^2 \vert x \vert^2 - R(x,I_\e(t))\\
&\geq (\overline K_1- 4 \overline M_1^2)\vert x \vert^2 \geq 0.\end{align*}
with $ b(0,x) =  \overline{L}_0 -\overline{M}_1 |x|^2 \geq  u_\e(0,x),$ $ x $ in $\R^d.$ Thus, $b$ is a super-solution of equation \fer{eq-u_e}, we conclude using the maximum principle.
 \end{itemize}
 \textbf{(ii) Let us prove the concavity estimate \fer{eq.D2u}.}\\\\
Let $\delta>0$ be small. We define $\phi_\delta$  $\in C^\infty(\R^+)$, as an increasing concave function, such that:
$$
\phi_\delta(r):=
\begin{cases}
r & \text{if } r\leq \delta^{-1}, \\[6pt]
\delta^{-1}+1 & \text{if } r\geq \delta^{-1}+2.
\end{cases}
$$
We define $\ued$ to be the solution of
\begin{align}
\label{ued}
\partial_t\ued-\e\Delta\ued&=\phi_\delta(\vert\nabla\ued\vert^2)+R(x,I_\e(t)),  \quad t>0, \ x \in \R^d,\\
\ued(0,x)&=u^0(x), \nonumber
\end{align}
and $\underline u$ to be the solution of
\begin{align}
\label{u_harm}
\partial_t\underline{u}-\e\Delta\underline{u}&= -2 \underline{K_1} \vert x \vert^2,  \quad t>0, \ x \in \R^d,\\
\underline u(0,x)&=u^0(x). \nonumber
\end{align}
\begin{itemize}[leftmargin=*] \item \textbf{We first show that $\underline u \leq u_{\e,\delta} \leq u_\e$}.\\
\textit{Let us show that $ \underline u \leq u_{\e,\delta}$}.\\
\\We use \fer{ued} to compute $$\partial_t \ued - \e \Delta \ued + 2 \underline K_1\vert x \vert^2 = 2 \underline K_1 \vert x \vert^2+ \phi_\delta(\vert \nabla \ued \vert ^2) + R(x,I_\e(t)).$$
Thanks to the definition of $\phi_\delta$, we know that $ \phi_\delta(r)\geq 0 $. According to the assumption \fer{asr} on $R$, we have $R(x,I_\e(t)) \geq -2\underline{K_1}\vert x \vert^2$.\\
We deduce that $\ued$ satisfies:
\begin{align*}
    \partial_t \ued - \e \Delta \ued + 2 \underline K_1 \vert x \vert^2 &\geq 0, \quad t>0, \; x \in \R^d, \\
    \ued(0,\cdot) = u^0 &\geq \underline{u}(0,\cdot).
\end{align*}
Thus, $\ued$ is a super-solution of the following PDE
\begin{align*}
\partial_t\underline{u}-\e\Delta\underline{u}=& -2 \underline{K}_1\vert x \vert^2,  \quad t>0, \ x \in \R^d,\\
\underline u(0,\cdot) = & \ u^0(\cdot).
\end{align*}
We define $v := \ued -\underline u $
and $v$ is a super-solution of:
\begin{align*}
    \partial_t v-\e\Delta v \geq& \ 0,  \quad t>0, \ x \in \R^d,\\
 v(0,\cdot)\geq &  \ 0.
\end{align*}
The application $v$ is a super-solution of a parabolic equation with bounded coefficients, and $0$ is the solution of the same equation, thus we obtain $\underline u \leq \ued$, using the comparison principle.\\\\
\textit{We now show that $ u_{\e,\delta} \leq u_\e$}.\\\\
We notice from the definition of $\phi_\delta$, that for all $r\geq0$, we have $\phi_\delta(r) \leq r$, using \fer{eq-u_e} we find
\begin{align*}
\partial_t u_\e-\e\Delta u_\e \geq&  \ \phi_\delta(\vert \nabla u_\e \vert^2) + R(x,I_\e(t)),  \quad t>0, \ x \in \R^d,\\
  u_\e(0,\cdot)\geq& \ \ued(0,\cdot).\\
\end{align*}
Thus, $u_\e$ is a super-solution of \fer{ued}. The application $\phi_\delta$ is bounded, we conclude using the maximum principle and we obtain $\underline u \leq \ued \leq u_\e.$
\end{itemize}
Next, we prove the following Lemma \fer{lemmeJM}.
\begin{lemma}
     \label{lemmeJM} Under assumptions \fer{aspsi}--\fer{asuIni},
     there exists constants $\overline C_{\e,\delta}, \underline{C}_{\e,\delta} >0$ such that
     $$-\underline{C}_{\e,\delta} \leq D^2 u_{\e,\delta} \leq \overline{C}_{\e,\delta}.$$
\end{lemma} 
\textbf{Proof.}
We want to show that $D^2 u_{\e,\delta}$ is bounded, but not necessarily uniformly with respect to $\e$ or $\delta.$
To this end, we define $v_{\e,\delta}$ and $w_{\e,\delta}$ as follow
\begin{align*}
\partial_t v_{\e,\delta} -\e \Delta v_{\e,\delta}& = R(x,I_\e (t)), \\
v_{\e,\delta}(0,\cdot) &= u_{\e}(0,\cdot),\\
    \partial_t w_{\e,\delta} -\e \Delta w_{\e,\delta} &= \phi_\delta(|\nabla u_{\e,\delta}|^2),\\
    w_{\e,\delta}(0,\cdot) &= 0,
\end{align*}
and we have \beq 
\label{somme} u_{\e,\delta} = v_{\e,\delta}+w_{\e,\delta},\eeq where $u_{\e,\delta}$ is the application defined in \fer{ued}.
\begin{itemize}[leftmargin=*]
    \item We can solve the PDE on $v_{\e,\delta}$ as the heat equation, and $D^2 v_\e$ satisfies 
    \begin{align*} D^2v_\e(t,x) =& \ \frac{1}{(4\pi t)^{\frac{d}{2}}} \int_{\R^d} e^{-(x-y)^2/(4t)} D^2 u_\e(0,y) \mathrm{d}y \\&+  \int_0^t \int_{\R^d} \frac{1}{(4\pi(t-s))^{\frac{d}{2}}} e^{-(x-y)^2/(4(t-s))} D^2 R(y,I_\e(s)) \mathrm{d}y \mathrm{d}s. \end{align*}
    We now use the assumption \fer{asuD2} on $u^0,$ to obtain that there exists $\underline C, \overline C,$ such that $$ -\underline C \leq \frac{1}{(4\pi t)^{\frac{d}{2}}} \int_{\R^d} e^{-(x-y)^2/(4t)} D^2 u_\e(0,y) \mathrm{d}y \leq - \overline{C}.$$
    We now estimate $J(t,x) :=  \int_0^t \int_{\R^d} \frac{1}{(4\pi(t-s))^{\frac{d}{2}}} e^{-(x-y)^2/(4(t-s))} D^2 R(y,I_\e(s)) \mathrm{d}y \mathrm{d}s. $ To this end, we use the change of variables $$ z = \frac{x-y}{2 \sqrt{t-s}},$$ and we find
    $$J(t,x) = \frac{1}{\pi^{\frac{d}{2}}}\int_0^t \int_{\R^d}  e^{-z^2} D^2 R(x-2z\sqrt{t-s},I_\e(s)) \mathrm{d}z \mathrm{d}s.$$
    Using the assumption \fer{asrD2} on $R$, we conclude that
$$ -\underline C(T) \leq D^2 v_{\e,\delta} \leq - \overline C(T),$$ Similarly, we show that there exists $\gamma(T)$ such that when $x$ tends to $+\infty$, \beq \label{equiv} |\nabla v_{\e,\delta} | \sim \gamma |x|.\eeq
\item The application $\phi_\delta(|\nabla u_{\e,\delta}|^2)$ is bounded, thus according to Ladyzhenskaya \cite{OL:88}, Chap.3, there exists $C(\e,\delta)>0$, such that $ |\nabla w_{\e,\delta} | \leq C(\e,\delta).$
\item We compute $\phi_\delta(|\nabla u_{\e,\delta}|^2)$ using \fer{somme}, and we find \begin{align*}
    \phi_{\delta}(|\nabla u_{\e,\delta}|^2) = \phi_{\delta}(|\nabla v_{\e,\delta}|^2+ 2\nabla v_{\e,\delta} \cdot \nabla w_{\e,\delta} +|\nabla w_{\e,\delta}|^2).
\end{align*}
Using the definition of $\phi_\delta$, the equivalence \fer{equiv} and the bound on $|\nabla w_{\e,\delta} |$, there exists $R_{\delta,\e}>0$, such as for all $|x| \geq R_{\delta,\e}, $ 
$$ \phi_{\delta} (|\nabla u_{\e,\delta}|^2) =1.$$
\item Thus, $w_{\e,\delta}$ satisfies an equation of the type:
$$\partial_t w_{\e,\delta} - \e \Delta w_{\e,\delta} =  f(t,x,\nabla w_{\e,\delta}),$$
with $f(t,x, \nabla w_{\e,\delta}) =1 $ if $|x| \geq R_{\delta,\e}$.
According to \cite{OL:88} Chap.4 and 5, since $v_{\e,\delta}$ is locally $C^{3+\alpha, 1 +\frac{\alpha}{2}}_{x,t}(\R^d\times\R^+)$, the map $(t,x) \mapsto f(t,x,p)$ is $C^{2+\alpha, 1+ \frac{\alpha}{2}}_{x,t}(\R^d\times\R^+)$ locally uniformly in $p$ (the Hölder constants do not depend on $p$). Thus $w_{\e,\delta}$ is $C^{2+ \alpha, 1+ \frac{\alpha}{2}}_{x,t}(\R^d\times\R^+)$ uniformly in $x.$ \end{itemize}
\begin{flushright}
    $ \square$
\end{flushright}
We continue the proof of Theorem \ref{t2.2}.
\begin{itemize}[leftmargin=*]
\item \textbf{Let us show that $-2\underline{M_1} \leq D^2\ued(t,x) \leq -2 \overline M_1$ : }\\\\
\textit{We prove the lower bound, the upper bound can be proven similarly.}
\\\\
We define $v := \partial_{\xi\xi}^2 \ued$, we differentiate \fer{ued} twice, to obtain
$$
\begin{array}{rll}
&\partial_t v-\e\Delta v-2\phi_\delta'(\vert\nabla\ued\vert^2)\nabla\ued.\nabla v-4\phi_\delta''(\vert\nabla\ued\vert^2)(\nabla\ued.\nabla (\partial_\xi \ued))^2\\=&2\phi_\delta'(\vert\nabla\ued\vert^2)\vert\nabla\partial_\xi\ued\vert^2+\partial_{\xi\xi}^2 R(x,I_\e(t)).
\end{array}
$$
We introduce $$
w(t,x):=\max_{\vert\xi\vert=1} \partial_{\xi\xi}^2 \ued(t,x).
$$ 
We prove that $w$ is a viscosity sub-solution of:
\beq 
\label{eqW} \partial_t w-\e\Delta w-2\phi_\delta'(\vert\nabla\ued\vert^2)\nabla\ued \cdot\nabla w=2\phi_\delta'(\vert\nabla\ued\vert^2)w^2-2\overline K_1.\eeq \\
Let $h$ be a test function such that $w-h$ admits a local maximum at $(t_0,x_0)$. Thanks to  Lemma \ref{lemme0}, there exists $\xi_0$ such that $w(t_0,x_0) = \partial_{\xi_0\xi_0}^2 u(t_0,x_0)$. The application $\partial_{\xi_0\xi_0}^2 u$ has a local maximum at $(t_0,x_0)$, because $w = \max_{\vert\xi\vert=1} \partial_{\xi\xi}^2 \ued$. Thus, we deduce that at $(t_0,x_0)$
$$\partial_t h - \e \Delta h -2\phi'_\delta(\vert \nabla \ued \vert^2) \nabla \ued \cdot \nabla h \leq \partial_t \partial_{\xi_0\xi_0}^2 u - \e \Delta \partial_{\xi_0\xi_0}^2 u -2\phi'_\delta(\vert \nabla \ued \vert^2) \nabla \ued \cdot \nabla \partial_{\xi_0\xi_0}^2 u
$$
$$ = 4\phi_\delta''(\vert\nabla\ued\vert^2)(\nabla\ued.\nabla (\partial_\xi \ued))^2 + 
2\phi_\delta'(\vert\nabla\ued\vert^2)\vert\nabla\partial_\xi\ued\vert^2+\partial_{\xi_0\xi_0}^2 R(x,I_\e(t)).$$\\
By definition, $\phi_\delta$ is concave, thus $\phi_\delta '' \leq 0$. According to the assumption \fer{asrD2} on $R$, we have $R_{\xi_0\xi_0} \leq -2\overline K_1$. Thanks to Lemma \ref{lemme0}, $\vert\nabla\partial_\xi\ued\vert^2 = w^2 $. Therefore, we obtain the following inequality
$$ \partial_t h - \e \Delta h -2\phi'_\delta(\vert \nabla \ued \vert^2) \nabla \ued \cdot \nabla h\leq 2\phi_\delta'(\vert\nabla\ued\vert^2)w^2-2\overline K_1 .$$
Moreover, $\overline w:=-2\overline M_1$ is a super-solution of the latter equation. We compute
$$2 \phi'_\delta (\vert \nabla \ued \vert ^2) \overline w ^2 -2 \overline K_1 = 8 \overline{M}_1^2 \phi'_\delta (\vert \nabla \ued \vert ^2) - 2 \overline K_1, $$
since $\phi_\delta$ is concave and  $\phi'_\delta(0) = 1$, we have $\phi'_\delta \leq 1$.\\
Thus, we deduce
$$2 \phi'_\delta (\vert \nabla \ued \vert ^2) \overline w ^2 -2 \overline K_1 \leq 2 (4 \overline{M}_1^2 - \overline{K}_1) \leq 0 .$$
Thanks to Lemma \ref{lemmeJM}, $w$ is bounded (not necessarly uniformly in $\e$ or $\delta$). Equation \fer{eqW} has bounded coefficients, thus using the maximum principle for viscosity solutions, we conclude that $D^2\ued \leq -2 \overline{M}_1.$ \\\\
\item \textbf{Let us prove that $(\ued)_\delta \xrightarrow[\delta \rightarrow 0]{} u_\e$ locally uniformly.}
The application $\ued$ is a viscosity solution of \fer{ued}, that we can write as $$\partial_t \ued + H_{\delta}(x,\ued,\nabla \ued, D^2 \ued) = 0$$
with $$H_{\delta}(x,u,p,M) = -\e tr(M) - \phi_\delta(\vert p \vert^2) -R(x,I_\e(t)).$$
-- $(\phi_\delta)_\delta$ converges locally uniformly to $Id$.\\
Thanks to the last points, \\
-- $ \ued$ is locally uniformly bounded in $\delta$, because $\underline{u} \leq \ued  \leq u_\e$.\\
-- $ \nabla \ued $ and $\partial_t \ued$ are locally uniformly bounded in $\delta$, because $-2 \underline{M_1} \leq D^2 \ued \leq -2 \overline{M_1}$ and $\underline{u} \leq \ued \leq u_\e$.\\
-- According to Arzela-Ascoli Theorem, there exists a subsequence $(\delta_k)_{k \in \N}$ such that \\ $u_{\e,\delta_k}\xrightarrow[k \rightarrow + \infty]{} v_\e$ locally uniformly.
Thanks to the stability Theorem of viscosity solutions of Hamilton-Jacobi equations (\cite{GB:94}, Theorem 2.1, p.21) $v_\e$ is a solution of \fer{eq-u_e}. We use the uniqueness of viscosity solution in $\R^d$ to obtain, $v_\e = u_\e$.\\
We let $\delta$  tend to $0$ and we obtain the wanted inequality  $-2 \underline{M}_1 \leq D^2 u_\e \leq -2 \overline{M}_1.$ \end{itemize}
\textbf{(iii) Let us prove that $ \|D^3 u_\e \|_{L^{\infty}([0,T]\times \R^d)} \leq C(T)$, \fer{nD3u}:}\\
We differentiate three times with respect to $x$ \fer{eq-u_e}, and obtain: 
\begin{align}
\label{u3}\partial_t \partial_{i,j,k}^3 u_\e - \e \Delta \partial_{i,j,k}^3 u_\e - 2 \nabla \partial_{i,j,k}^3 u_\e \cdot \nabla  u_\e = \ &
  2 \nabla \partial_{j,k}^2 u_\e \cdot \nabla \partial_i u_\e + 2 \nabla \partial_{j} u_\e \cdot \nabla \partial_{i,k}^2 u_\e \\ & + 2 \nabla \partial_{i,j}^2 u_\e \cdot \nabla \partial_{k} u_\e + \partial_{i,j,k}^3 R. \nonumber \end{align}
The function $\partial_{i,j,k}^3 u$ satisfies a parabolic system of equations \fer{u3}, the terms of the right-hand-side have bounded coefficients, and $D^3 u_0$ is uniformly bounded with respect to $\e$, thus according to the maximum principle, we have: $$ \|D^3 u_\e \|_{L^{\infty}([0,T]\times \R^d)} \leq C(T).$$
\textbf{(iv) One can prove the estimates for the derivatives of order 4, 5, and 6 similarly.}
\begin{flushright}
    $\square$
\end{flushright}
\section{\texorpdfstring{An approximation of $I_\e$}{}}
\label{s3}
In this section, we give an approximation of $I_\e$ which leads to an estimate of $R(x_\e(t),I_\e(t))$, where $x_\e(t)$  defined in \eqref{xe}, is the maximum of $u_\e(t,\cdot)$. We will use this estimate in the next section to prove the Proposition \ref{item1}. To this end, we introduce the set $\Omega := \{x \in \R^d, \ R(x,0) >0\}$ and we define $\mathcal I :\Omega \to \R^+$ as follows. Let $y\in \Omega$, then $\mathcal I (y)\in \R^{+*}$ is the unique nonnegative constant such that  $ R(y,\mathcal I (y))=0$. The function $\mathcal I$ is well defined on $\Omega$ because $R$ is strictly decreasing and tends to $-\infty$ as $I$ approaches $+\infty$.
 \begin{lemma}
\label{lemme1}
Under assumptions \fer{aspsi}--\fer{asuIni}, there exists a constant $C(T) = C>0$, independant of $\e$, such that for $t\in [0,T],$ we have
 $$
 |I_\e(t) - \mathcal I \big(x_\e(t)\big)| \leq  C\e,
 $$
 and consequently
 $$
 |R\big(x_\e(t),I_\e(t)\big) | \leq C\e.
 $$
 \end{lemma}

 \proof
 We give the proof for $d=1$. The argument can be adapted easily to the general case with $d\geq 1$. Let fix $T>0.$
 \\
 
 \noindent
 \textbf{(i)} We first need to check that for all $t \in [0,T],$ we have $x_\e(t) \in \Omega.$ Let $t$ be in $[0,T].$ It has been shown previously in \cite{GB.SM.BP:09}, \cite{AL.SM.BP:10} that $u_\e$ converges locally uniformly to $u$. By definition of $x_\e(t)$ and $\overline x(t)$, they respectively maximize $u_\e(t,\cdot)$ and $u (t,\cdot)$. Thus, $x_\e(t)$ converges to $ \bar x(t)$ along subsequences. We know by definition of $\mathcal I$ and by \cite{SM.JR:16} that
 $$R\big(\overline{x}(t),\mathcal{I}(\overline{x}(t))\big) = R\big(\overline{x}(t),I(t)\big) = 0.$$
 Thanks to \fer{asrDi}, we observe that $I(t) = \mathcal{I}\big(\bar x(t)\big)$.\\
We know that $I(t)>0$, thus because $R$ is strictly decreasing, we have $$R\Big(\overline{x}(t),\frac{I(t)}{2}\Big) >0.$$
Since $u_\e$ converges locally uniformly to $u$, for $\e$ sufficiently small, we also have
$$R\Big(x_\e(t),\frac{I(t)}{2}\Big) >0. $$
We conclude using \fer{asrDi} and $I(t)>0$ again, that $$R\big(x_\e(t),0\big) >0, \text{ \textit{ie.} $x_\e(t) \in \Omega$}.$$
\\\textbf{(ii)} We differentiate with respect to time the equation \fer{paraI}. Then we use the expression of $\e\partial_t n_\e$ from  equation \fer{para} and integrate by parts to obtain
$$\e\frac{\mathrm{d} }{\mathrm{d}t}I_\e(t) = \e^2 \int_{\R} \Delta n_\e(t,x) \psi(x) \mathrm{d}x + \int_{\R} n_\e(t,x)R \big(x,I_\e(t)\big)\psi(x)\mathrm{d}x.$$
We define $J_\e:= \log(I_\e).$ Using the Hof-Cole transformation \fer{Hopf}, we have
\begin{align}
\label{ODE-J_e}
\e \frac{\mathrm{d} }{\mathrm{d}t} J_\e(t) =& \ \e^2 \frac{\int_{\R} n_\e(t,x) \Delta \psi (x)\mathrm{d}x}{\int_{\R} n_\e(t,x) \psi (x)\mathrm{d}x} + R \big(x_\e(t),I_\e(t) \big) \\& + \frac{\int_{\R}
e^{\frac{u_\e(t,x)-u_\e(t,x_\e(t))}{\e}}\Big(R \big(x,I_\e(t) \big)- R \big(x_\e(t),I_\e(t) \big)\Big) \psi(x)\mathrm{d}x}{\int_{\R} e^{\frac{u_\e(t,x)-u_\e(t,x_\e(t))}{\e}} \psi (x)\mathrm{d}x} \nonumber.\end{align}
Using \fer{aspsi}, we control the first term of \fer{ODE-J_e} and obtain \beq \label{est1}\e^2 \frac{\int_{\R} n_\e(t,x) \Delta \psi (x)\mathrm{d}x}{\int_{\R} n_\e(t,x) \psi (x)\mathrm{d}x} = O(\e^2).\eeq
We then prove that the last term of the right-hand side of \fer{ODE-J_e} is of order $\e$. Firstly, we will find an equivalent of the denominator and secondly show that the numerator is of order $\e.$\\\\
Let us show that \beq \label{laplace} \int_{\R} e^{\frac{u_\e(t,x)-u_\e(t,x_\e(t))}{\e}} \psi (x) \mathrm{d}x =  \sqrt{\frac{2\pi \e}{\vert D^2 u(t,\overline{x}(t))\vert }} \psi \big(\overline{x}(t)\big) + o(\sqrt \e). \eeq
We use the fact that $x_\e(t)$ maximizes $u_\e(t,\cdot),$ and obtain the following Taylor-Lagrange expansion
$$\forall x \in \R, \ \exists z \in e(x,x_\e(t)), \ u_\e(t,x) = u_\e(t,x_\e(t)) + \f 12 D^2u_\e(t,x_\e(t))(x-x_\e(t))^2 + \f 16 D^3u_\e(t,z)(x-x_\e(t))^3.$$
We introduce $\alpha >0,$ to be chosen later and we split the integral in two
\begin{align*}
    A_1 &:= \int_{\vert x-x_\e(t) \vert  \leq \e^\alpha } e^{\frac{u_\e(t,x)-u_\e(t,x_\e(t))}{\e}} \psi (x)\mathrm{d}x,\\
    A_2  &:= \int_{\vert x-x_\e(t) \vert > \e^\alpha } e^{\frac{u_\e(t,x)-u_\e(t,x_\e(t))}{\e}} \psi (x)\mathrm{d}x.
\end{align*}
To estimate $A_1$, we make the change of variables $y = \frac{x-x_\e(t)}{\sqrt{\e}}$ and obtain
\begin{align*}
     A_1 & =  \sqrt \e  \int_{|y|\leq \e^{\alpha-\f 12}} e^{\frac{1}{2} D^2u_\e(t,x_\e(t))y^2 + \f{\sqrt{\e}}{6} D^3u_\e(t,z)y^3 }\psi(x_\e(t)+\sqrt \e y)\mathrm{d}y.
\end{align*}
We want to write an asymptotic expansion of the exponential. To this end, we use that $D^3u_\e$ is uniformly bounded thanks to \fer{nD3u}, and choose $\alpha$ such that $\sqrt \e y^3 = O(\e^{3\alpha -1})$ is small, ie. $\alpha>\frac{1}{3}$ and compute
\begin{align*}
    A_1  = & \ \sqrt \e  \int_{|y|\leq \e^{\alpha-\f 12}} e^{\frac{1}{2} D^2u_\e(t,x_\e(t))y^2} \Big( 1+ \frac{\sqrt \e}{6} D^3u_\e(t,z)y^3 + O(\e y^6)\Big)\\ & \times \Big( \psi \big(x_\e(t)\big) + \sqrt{\e} y \psi' \big(x_\e(t) \big) + O(y^2\e)\Big) \mathrm{d}y \\
      =& \ \sqrt \e \psi\big(x_\e(t)\big)\int_{|y|\leq \e^{\alpha-\f 12}}e^{\frac{1}{2} D^2u_\e(t,x_\e(t))y^2}\mathrm{d}y + \frac{\e}{6}\psi\big(x_\e(t)\big) \int_{|y|\leq \e^{\alpha-\f 12}}e^{\frac{1}{2} D^2u_\e(t,x_\e(t))y^2} D^3 u_\e (t,z) y^3\mathrm{d}y\\
       & +\sqrt \e \psi\big(x_\e(t) \big)\int_{|y|\leq \e^{\alpha-\f 12}}e^{\frac{1}{2}D^2u_\e(t,x_\e(t))y^2}O(\e y^6) \mathrm{d}y+ \mathcal{E}(\e,\alpha),
\end{align*}
 where $\mathcal{E}(\e,\alpha)$ is the rest of the expansion. 
We will compute and estimate the first terms of the latter equation, before we prove that $\mathcal{E}(\e,\alpha) = o(\sqrt \e).$\\
We choose $\frac{1}{3}<
\alpha < \frac{1}{2}, $  such that $\lim_{\e \to 0} \e^{\alpha-\frac{1}{2}} = +\infty,$ and thanks to \fer{eq.D2u} we have:  $$ \sqrt \e \psi \big(x_\e(t) \big) \int_{|y|\leq \e^{\alpha-\f 12}}  e^{\frac{1}{2} D^2 u_\e(t,x_\e(t))y^2 } \mathrm{d}y = \psi \big(\overline{x}(t) \big)\sqrt{\frac{2\pi\e}{\vert D^2 u(t,\overline{x}(t))\vert }} + o(\sqrt \e).$$
The second term of $A_1$ is of order $o(\sqrt \e)$. Indeed, for $\alpha >\frac{1}{3}$, $|\e y^3| \leq \e^{3\alpha -\frac{1}{2}} = o(\sqrt \e).$ We then use  Theorem \ref{t2.2}, which gives us that $D^3 u_\e$  is uniformly bounded \fer{nD3u}, and the concavity estimate on $u_\e$ \fer{eq.D2u} to obtain $\lim_{\e \to 0} \int_{0}^{\e^{\alpha-\f 12}} e^{\frac{1}{2} D^2 u_\e(t,x_\e(t))y^2 } \mathrm{d}y \leq \int_{0}^{+\infty} e^{-\overline M_1 y^2} \mathrm{d}y<+\infty.$
\\Similarly, the third term of the right-hand side is $o(\sqrt \e)$. We now want to make sure that the last term is small enough, we write $\mathcal{E}(\e,\alpha)$
\begin{align*}
    \mathcal{E}(\e,\alpha) =& \ \e \psi' \big(x_\e(t) \big)\int_{|y|\leq \e^{\alpha-\f 12}}  e^{\frac{1}{2} D^2 u_\e(t,x_\e(t))y^2 } y \Big(1+\frac{\sqrt \e}{6}D^3 u_\e (t,z) y^3 + O(\e y^6)\Big)\mathrm{d}y \\ 
    &+ \e \sqrt{\e}\int_{|y|\leq \e^{\alpha-\f 12}}  e^{\frac{1}{2} D^2 u_\e(t,x_\e(t))y^2 } O(y^2)\Big(1+ \frac{\sqrt \e}{6} D^3 u_\e (t,z)y^3 + O(\e y^6)\Big) \mathrm{d}y.
\end{align*}
Using the fact that $\psi\in W^{2,\infty}(\R^d)$ thanks to \fer{aspsi}, as well as Theorem \ref{t2.2} which gives us the estimate $D^2 u_\e \leq -2 \bar L_1$ \fer{eq.D2u}, and the uniform bound on $D^3 u_\e$ \fer{nD3u} and recalling that we chose $\frac{1}{3} < \alpha < \frac{1}{2}$, we deduce that $\mathcal{E}(\e,\alpha) = o(\sqrt{\e}).$\\\\
We now estimate $A_2$.
One can make the same change of variables $y = \frac{x-x_{\e}(t)}{\sqrt{\e}}$ to obtain
$$ A_2 = \sqrt \e  \int_{|y|> \e^{\alpha-\f 12}} e^{\frac{1}{2} D^2u_\e(t,x_\e(t))y^2 + \f{\sqrt{\e}}{6} D^3(t,z)y^3 }\psi \big(\sqrt \e y + x_\e(t)\big)\mathrm{d}y.$$
Using again \fer{aspsi} \fer{eq.D2u} \fer{nD3u}, we obtain that the application under the integral is integrable. We also have $\e^{\alpha-\f 12} 
  $ goes to $+\infty$ when $\e$ approaches 0. Thus, we have $A_2 = o(\sqrt \e).$
Finally, we obtain
$$\int_{\R} e^{\frac{u_\e(t,x)-u_\e(t,x_\e(t))}{\e}} \psi (x)\mathrm{d}x = A_1 + A_2= \sqrt{\frac{2\pi \e}{\vert D^2 u(t,\overline{x}(t))\vert }} \psi \big(\overline{x}(t) \big) + o(\sqrt \e). $$
We now prove that the numerator of the last term in \fer{ODE-J_e}, is of order $\e \sqrt \e$, ie.
\beq \label{est3} \Big|\int_{\R} e^{\frac{u_\e(t,x)-u_\e(t,x_\e(t))}{\e}}\Big(R(x,I_\e(t))- R\big(x_\e(t),I_\e(t)\big)\Big) \psi(x)\mathrm{d}x \Big| \leq C\e \sqrt \e.\eeq We proceed as before and split the integral in two, with $\alpha$ to be chosen later
\begin{align*}   
    &A_1 :=\int_{|x-x_\e|\leq  \e^{\alpha} } e^{\f {u_\e(t,x)-u_\e(t,x_\e(t))}{\e}} \psi(x) \Big( R\big(x,I_\e(t)\big) - R\big(x_\e(t), I_\e(t)\big) \Big)\mathrm{d}x,\\
    &A_2 := \int_{|x-x_\e|> \e^{\alpha} }  e^{\f {u_\e(t,x)-u_\e(t,x_\e(t))}{\e}} \psi(x) \Big( R\big(x,I_\e(t)\big) - R\big(x_\e(t), I_\e(t)\big) \Big)\mathrm{d}x.
\end{align*}
We write the following Taylor-Lagrange expansions:
\begin{align*}
    &\forall x \in \R, \ \exists z \in e(x,x_\e(t)), \ u_\e(t,x) = u_\e (t,x_\e(t)) + \f 12 D^2u_\e(t,z)(x-x_\e(t))^2,  \\
    &\forall x \in \R, \ \exists z' \in e(x,x_\e(t)), \ R \big(x,I_\e(t) \big)  =  R\big(x_\e(t),I_\e(t)\big) + \nabla_x R\big(x_\e(t),I_\e(t)\big)\cdot (x-x_\e(t)) \\
    & \quad \quad \quad \quad \quad \quad \quad \quad \quad \quad \quad \quad \quad \quad \quad \quad \quad  \ + \f12 D^2 R \big(z', I_\e(t)\big) (x-x_\e(t))^2.
\end{align*}
We start by controlling $A_2$. Using the assumption  \fer{aspsi} on $\psi$, the uniform boundedness of $x_\e(t)$, the assumptions \eqref{asr} \fer{asrD2} on $R$, and the estimate \fer{eq.D2u} on $D^2 u_\e$,  we obtain the following inequalities
\begin{align*}
    \vert A_2 \vert  \leq& \ \Big\vert \int_{|x-x_\e|> \e^{\alpha}} e^{\frac{D^2u_\e(t,z)(x-x_\e(t))^2}{2\e}} \psi(x) \nabla_x R \big(x_\e(t),I_\e(t)\big)(x-x_\e(t))\mathrm{d}x \Big \vert  \\
    & + \f 12 \Big \vert \int_{|x-x_\e|> \e^{\alpha}} e^{\frac{D^2u_\e(t,z)(x-x_\e(t))^2}{2\e}} \psi(x) D^2 R\big(z',I_\e(t)\big)(x-x_\e(t))^2 \mathrm{d}x \Big \vert \\
    \leq & \ C \| \psi \|_\infty \int_{|x-x_\e|> \e^{\alpha}} \exp \Big(\frac{-\overline{M}_1\vert x-x_\e(t) \vert^2}{\e} \Big) \vert x-x_\e(t) \vert \mathrm{d}x \\
    & + C \| \psi \|_\infty \| D^2 R \|_\infty \int_{|x-x_\e|> \e^{\alpha}} \exp \Big(\frac{-\overline{M}_1 \vert x-x_\e(t) \vert^2}{\e} \Big) \vert x-x_\e(t) \vert ^2 \mathrm{d}x \\
  \leq & \ 2 C \| \psi \|_\infty \int_{\e^{\alpha}}^{+\infty} x e^{-\frac{\overline M_1 x^2}{\e}} \mathrm{d}x + 2 C \| \psi \|_\infty \| D^2 R \|_\infty \int_{\e^\alpha}^{+ \infty} x^2 e^{-\frac{\overline M_1 x^2}{\e}} \mathrm{d}x.
\end{align*}
We can compute the two integrals of the right-hand side explicitly and we obtain
\begin{align*}
    &\int_{\e^{\alpha}}^{+\infty} x e^{-\frac{\overline M_1 x^2}{\e}} dx = \frac{ \e}{2\overline{M}_1} e^{-\overline{M}_1\e^{2\alpha -1}},\\
&\int_{\e^\alpha}^{+ \infty} x^2 e^{-\frac{\overline M_1 x^2}{\e}} \mathrm{d}x =  \frac{\e^{\alpha+1}}{2\overline{M}_1} e^{- \overline M_1 \e^{2\alpha-1}} + \frac{\e}{2\overline{M}_1}\int_{\e^\alpha}^{+ \infty} e^{-\frac{\overline M_1 x^2}{\e}} \mathrm{d}x \leq \frac{\e^{\alpha+1}}{2\overline{M}_1} e^{-  \overline M_1 \e^{2\alpha-1}} +  \frac{\e}{4\overline M_1} \sqrt{\frac{\pi \e}{\overline{M}_1}}.
\end{align*}
Plugging the results of the integrals in the latter inequality, we obtain
\begin{align*}
    \vert A_2 \vert & \leq C  \e e^{-  \overline M_1 \e^{2\alpha-1}} + C \e^{\alpha+1} e^{- \overline M_1 \e^{2\alpha-1}} +  C \e \sqrt{\e}.
\end{align*}
Thus, we have $ A_2 = O(\e \sqrt \e ), \text{ with $0<\alpha < \f 12$}.$\\\\
To estimate $A_1$, we write a Taylor-Lagrange expansion to the next order. There exists $z', z'' \in e(x,x_\e(t))$, such that
\begin{align*}
    A_1  =& \ \int_{|x-x_\e|\leq \e^{\alpha}} e^{\f{\f12 D^2 u_\e(t,x_\e(t))(x-x_\e(t))^2 +\f16 D^3 u_\e(t,z'')(x-x_\e(t))^3 }{\e}}  \nabla_x R\big(x_\e(t),I_\e(t)\big)\cdot (x-x_\e(t)) \mathrm{d}x \\
    & + \int_{|x-x_\e|\leq \e^{\alpha}} e^{\f{\f12 D^2 u_\e(t,x_\e(t))(x-x_\e(t))^2 +\f16 D^3 u_\e(t,z'')(x-x_\e(t))^3 }{\e}} \f 12D^2 R\big(z',I_\e(t)\big) (x-x_\e(t))^2 \mathrm{d}x.
\end{align*}
After the change of variables $y = \frac{x-x_\e(t)}{\sqrt \e}$, and for $ \frac{1}{3} < \alpha < \frac{1}{2}$, we have
\begin{align*}
    A_1 =& \ \e \int_{|y|\leq \e^{\alpha-\f12}} e^{\f12 D^2 u_\e(t,x_\e(t))y^2 +\frac{\sqrt{\e}}{6}  D^3 u_\e(t,z'')y^3 } \left[ \nabla_x R\big(x_\e(t),I_\e(t)\big)\cdot y
 +\frac{\sqrt \e}{2} D^2 R\big(z',I_\e(t)\big) y^2\right]\\
 &\times \psi(x_\e(t) + \sqrt \e y)\mathrm{d}y \\
 =& \ \e \int_{|y|\leq \e^{\alpha-\f12}}e^{\f12 D^2 u_\e(t,x_\e(t))y^2}\big(1+ \frac{\sqrt \e}{6} D^3 u_\e(t,z'')y^3 + O(\e y^6)\big)\\ & \times\big[ \nabla_x R\big(x_\e(t),I_\e(t)\big) \cdot y + \frac{\sqrt \e}{2}  D^2 R\big(z',I_\e(t)\big) y^2 \big]
  \big( \psi(x_\e(t)) + \sqrt \e y \psi '(x_\e(t)) +O(\e y^2)\big) \mathrm{d}y\\
      =& \ \e \psi \big(x_\e(t) \big) \nabla_x R\big(x_\e(t),I_\e(t)\big)\int_{|y|\leq \e^{\alpha-\f12}} e^{\f12 D^2 u_\e(t,x_\e(t))y^2} y \mathrm{d}y 
 \\& + \e \psi\big(x_\e(t) \big) \nabla_x R\big(x_\e(t),I_\e(t)\big)\int_{|y|\leq \e^{\alpha-\f12}} e^{\f12 D^2 u_\e(t,x_\e(t))y^2} O(\e |y|^7) \mathrm{d}y
 \\& + \frac{\e \sqrt \e}{2} \psi\big(x_\e(t) \big) \int_{|y|\leq \e^{\alpha-\f12}} e^{\f12 D^2 u_\e(t,x_\e(t))y^2} D^2 R\big(z',I_\e(t)\big) y^2 \mathrm{d}y\\
 &+\frac{\e \sqrt \e}{6} \psi\big (x_\e(t) \big) \nabla_x R\big(x_\e(t),I_\e(t)\big) \int_{|y|\leq \e^{\alpha-\f12}} e^{\f12 D^2 u_\e(t,x_\e(t))y^2} D^3 u_\e \big(z'',I_\e(t)\big) y^4 \mathrm{d}y \\
 &+ \e \sqrt{\e} \psi'\big(x_\e(t) \big) \nabla_x R\big(x_\e(t),I_\e(t)\big)  \int_{|y|\leq \e^{\alpha-\f12}} e^{\f12 D^2 u_\e(t,x_\e(t))y^2} y^2 \mathrm{d}y + \mathcal{E}(\e,\alpha).
\end{align*}
The first term of the right-hand side is $0$ because the interval of integration is symmetric and the function is odd. The second term, $\e \psi\big(x_\e(t)\big) \nabla_x R\big(x_\e(t),I_\e(t)\big)\int_{|y|\leq \e^{\alpha-\f12}} e^{\f12 D^2 u_\e(t,x_\e(t))y^2} O(\e |y|^7) \mathrm{d}y$ is of order $O(\e ^2),$ indeed the integral $\int_{\R} e^{\f12 D^2 u_\e(t,x_\e(t))y^2}|y|^7 \mathrm{d}y$ is finite, thanks to the concavity estimate \fer{eq.D2u} on $u_\e$. Using the concavity estimate \fer{eq.D2u} on $u_\e$, the assumptions \fer{asrD2} on $D^2R$ and $D^3 R$, we conclude that the other terms are of order $O(\e \sqrt \e).$
Hence, $$\Big |\int_{\R} e^{\frac{u_\e(t,x)-u_\e(t,x_\e(t))}{\e}}\Big(R\big(x,I_\e(t)\big)- R\big(x_\e(t),I_\e(t)\big)\Big) \psi(x)\mathrm{d}x  \Big | \leq C\e \sqrt \e.$$ Finally, we combine the last three estimates \fer{est1} \fer{laplace} \fer{est3} with  the equation \fer{ODE-J_e}, and we obtain the following ODE: \beq\label{ODE-Je}\e \f{\mathrm{d} }{\mathrm{d} t }J_\e(t)= R\big(x_\e(t),I_\e(t)\big)+O(\e).\eeq
\textbf{(iii)} We denote $\mathcal J := \log(\mathcal I)$. By definition of $\mathcal{I}$, we have \beq \label{Re=0}R\big(x_\e(t),\mathcal{I}(x_\e(t))\big) = 0.\eeq
We differentiate the latter equation \fer{Re=0} with respect to time and obtain 
$$ \nabla_x R\Big(x_\e(t), \mathcal{I}\big(x_\e(t)\big)\Big) \dot{x}_\e(t) + \frac{\partial}{\partial I}R\big(x_\e(t), \mathcal{I}(x_\e(t))\big) \frac{\mathrm{d}}{\mathrm{d}t}\mathcal{I}\big(x_\e(t)\big) = 0,$$
thus, $$\frac{\mathrm{d}}{\mathrm{d}t} \mathcal J\big(x_\e(t)\big) = -\frac{ \nabla_x R\Big(x_\e(t), \mathcal{I}\big(x_\e(t)\big)\Big) \dot{x}_\e(t)}{\mathcal{I}\big(x_\e(t)\big)\frac{\partial}{\partial I}R\Big(x_\e(t), \mathcal{I}\big(x_\e(t)\big)\Big)}.$$
We know that $x_\e$ converges to $\overline x$ and that $ 0<I_m\leq \mathcal{I}(\overline x(t)) = I(t) \leq 2I_M$. We deduce that for $\e$ sufficiently small, we have $0<\frac{I_m}{2} \leq \mathcal{I}(x_\e(t)) \leq 2 I_M$ for all $t$ in $[0,T]$. Now, we use assumption \fer{asrDi}, and obtain $$-\frac{ \nabla_x R\Big(x_\e(t), \mathcal{I}\big(x_\e(t)\big)\Big)}{\mathcal{I}\big(x_\e(t)\big)\frac{\partial}{\partial I}R\Big(x_\e(t), \mathcal{I}\big(x_\e(t)\big)\Big)} = O(1).$$
We need to verify that $\dot{x}_\e(t)=O(1).$ To this end, we use that by definition $x_\e(t)$ maximizes $u_\e(t,\cdot)$, which leads to
$$\nabla u_\e(t,x_\e(t)) = 0.$$
We differentiate with respect to time to obtain 
    $$\frac{\mathrm{d}}{\mathrm{d}t}\nabla u_\e(t,x_\e(t)) =\nabla \partial_t  u_\e(t,x_\e(t)) + D^2 u_\e(t,x_\e(t))\dot x_\e(t) = 0$$
Thanks to \fer{para} we have for all $i\in \{1,...,d\},$ $$\partial_i \partial_t u_\e = \e\partial_i \Delta u_\e + 2 \partial_i(\vert \nabla u_\e \vert) \vert \nabla u_\e \vert + \partial_i R(x,I_\e(t)). $$
Hence, we have \beq \label{dx_e} \dot{ x}_\e(t) = \left( -D^2u_\e \big(t,  x_\e(t)\big) \right)^{-1}  \left[ \nabla_x R\big( x_\e(t),I_\e(t)) \big) +\e \nabla \Delta u_\e \right]. \eeq
By Theorem \ref{thm:approx} and the uniform boundedness of $x_\e(t)$, we have $$ \dot{x}_\e(t) = O(1).$$
Finally, for all $t \in [0,T]$, we have $$ \e\frac{\mathrm{d}}{\mathrm{d}t} \mathcal J\big(x_\e(t)\big) = R\Big(x_\e(t), \mathcal{I}\big(x_\e(t)\big)\Big) + O(\e).$$
We subtract the latter equation from the equation \fer{ODE-Je} and we obtain for all $t \in [0,T],$ $$ \e\frac{\mathrm{d}}{\mathrm{d}t} \Big( J_\e(t)- \mathcal J\big(x_\e(t)\big) \Big) = R\big(x_\e(t),I_\e(t)\big) - R\big(x_\e(t), \mathcal{I}\big(x_\e(t)\big)\big) + O(\e).$$
We multiply by $ \mathrm{sgn}\big( J_\e(t)- \mathcal J\big(x_\e(t)\big)\big) $ and use the fact that $R$ is decreasing with respect to $I$, thanks to the assumption \fer{asrD2}, to obtain the inequality 
$$ \e \frac{\mathrm{d}}{\mathrm{d}t} \vert J_\e(t)- \mathcal J\big(x_\e(t)\big) \vert \leq -\vert R\big(x_\e(t),I_\e(t)\big) - R\big(x_\e(t), \mathcal{I}\big(x_\e(t)\big)\big) \vert  + O(\e).$$
We denote $\tilde{R}(x,J) := R(x, e^J) $, and we use the mean value theorem to obtain
$$\e \frac{\mathrm{d}}{\mathrm{d}t} \vert J_\e(t)- \mathcal J\big(x_\e(t)\big) \vert \leq \min_{x,J} \frac{\partial }{\partial J}\tilde R(x,J)  \vert J_\e(t)- \mathcal J\big(x_\e(t)\big) \vert  + O(\e).$$
We apply Grönwall's lemma to the last inequality 
$$ \vert J_\e(t)- \mathcal J\big(x_\e(t)\big) \vert \leq  A e^{-\frac{Ct}{\e}} + O(\e).$$
Because $I_\e$ and $\mathcal{I}(x_\e(\cdot))$ are bounded in $[0,T]$, we have $$ \vert I_\e(t)- \mathcal I\big(x_\e(t)\big) \vert \leq  \tilde{A} e^{-\frac{Ct}{\e}} + O(\e),$$
i.e. for $T >0$ fixed, for $ t$ in $ [0,T],$ $ |I_\e(t) - \mathcal I\big(x_\e(t)\big) |\leq  C\e,$
and consequently, for $ t$ in $ [0,T],$ $R\big(x_\e(t),I_\e(t)\big) = R\big(x_\e(t), \mathcal{I}\big(x_\e(t)\big)\big)+O(\e) = O(\e).$ This completes the proof of  Lemma \ref{lemme1}.
\begin{flushright} $\square$ \end{flushright}
\section{\texorpdfstring{Error estimate to the order $\e.$}{}}
  In this section, we will prove the Proposition \ref{item1} using Lemma \ref{lemme1}. This proposition is a first step before proving the Theorem \ref{thm:approx}. To this end, we adapt some methods used in \cite{SM.JR:16}.
  \begin{proposition}
      \label{item1}
      We assume \fer{aspsi}--\fer{asuIni}.
Let $n_\e$ be the solution of \fer{para} and $u_\e$ be defined by \fer{Hopf}. \\
There exists a constant $C(T) = C>0$, independant of $\e$, such that we have the following error estimate
\beq
\label{errorI}
\|I_\e  - I\|_{L^\infty([0,T])}  \leq C \e, \eeq
\beq \label{errorx}\|x_\e - \overline x \|_{L^\infty([0,T])} \leq C \e, \eeq 
\beq
\label{erroru}
\| u_\e - u- \e \log (\f{r}{ \e^{\f d 2}})\|_{L^\infty_t W^{4,\infty}_x([0,T]\times\R^d)}\leq  C\e,
\eeq
where $x_\e(t)$ is the maximum point of $u_\e$:
 $$
 \max_{x\in \R^d}u_\e(t,x) =u_\e(t,x_\e(t)).
 $$
  \end{proposition}
 \subsection{\texorpdfstring{System of equations satisfied by $(u_\e,I_\e,x_\e)$}{}}
 \noindent
 We first notice that $(u_\e,I_\e,x_\e)$ solves the following system
 \beq
\label{system-e}
\begin{cases}
R \big( x_\e(t),I_\e(t)) \big)=O(\e),& \text{for $t\in [0,T]$},\\
\dot{ x}_\e(t) = \big( -D^2u_\e \big(t,  x_\e(t)\big) \big)^{-1}  \left[ \nabla_x R\big( x_\e(t),I_\e(t)) \big) +\e \nabla \Delta u_\e \right],& \text{for $t\in [0,T]$},\\
\p_t u_\e = | \nabla  u_\e |^2 +R(x,I_\e) + \e \Delta u_\e,& \text{in $ [0,T] \times \R^d$},
\end{cases}
\eeq
with initial conditions
$$
u_\e(0,x)=u^0+  \e \log \Big(\f{r}{ \e^{\f d 2}}\Big) ,\quad  I_\e(0)=I_0+\e J_0+O(\e^2), \quad  x_\e(0)=x_0+O(\e).
$$
Indeed, we proved in Lemma \ref{lemme1} that $R\big(x_\e(t),I_\e(t)\big)=O(\e)$. The second line \fer{dx_e}  was proved in the last section.
 This system of equation \fer{system-e} allows us to compare $(u_\e,I_\e,x_\e)$ with $(u,I,x)$ which solves the system
 \beq
\label{system-u}
\begin{cases}
R \left( \bar x(t),I(t)) \right)=0,& \text{for $t\in [0,T]$},\\
\dot{ \bar x}(t) = \left( -D^2u \big(t,  \bar x(t)\big) \right)^{-1}   \nabla_x R\big( \bar x(t),I(t)) \big),& \text{for $t\in [0,T]$},\\
\p_t u = | \nabla  u |^2 +R(x,I), & \text{in $ [0,T] \times \R^d$}.
\end{cases}
\eeq 
To this end, we define
 $$
 v_\e(t,x)= u_\e(t,x)+\e \log (\f{r}{ \e^{\f d 2}}).
 $$
In this way, we remove the term of order $\e \log \e$, so that $(v_\e,I_\e,x_\e)$ solves
\beq
\label{system-em}
\begin{cases}
R \left( x_\e(t),I_\e(t)) \right)=O(\e),& \text{for $t\in [0,T]$},\\
\dot{ x}_\e(t) = \left( -D^2v_\e \big(t,  x_\e(t)\big) \right)^{-1}  \left[ \nabla_x R\big( x_\e(t),I_\e(t)) \big) +\e \nabla \Delta v_\e \right],& \text{for $t\in [0,T]$},\\
\p_t v_\e = | \nabla  v_\e |^2 +R(x,I_\e) + \e \Delta v_\e,& \text{in $ [0,T] \times \R^d$}.
\end{cases}
\eeq
with initial conditions
$$
v_\e(0,x)=u^0,  \quad  I_\e(0)=I_0+\e J_0+O(\e^2), \quad  x_\e(0)=x_0+O(\e).
$$
\subsection{\texorpdfstring{An approximation for $x_\e$, $I_\e$ and $u_\e$, of order $O(\e)$}{}}
In this subsection, we prove Proposition \ref{item1}.\\\\
\textbf{Proof.}
We fix $T>0$, and $t\in [0,T].$\\\\
\textbf{(i)} We start by proving the following inequality $$\vert I_\e(t) - \mathcal I(x_\e(t))\vert \leq C \vert x_\e(t) - \overline x(t)\vert + C\e.$$
By definition of $\mathcal{I}$, we have for all $y \in \Omega$, $ R\big(y, \mathcal I(y)\big) = 0,$
thus, $$R \big(\overline{x}(t), \mathcal I(\overline{x}(t))\big) - R\big(\overline{x}(t), \mathcal I(x_\e(t))\big)= R\big(x_\e(t), \mathcal I(x_\e(t))\big)-R\big(\overline{x}(t), \mathcal I(x_\e(t))\big).$$
Thanks to Lemma \ref{lemme1} we have $R\big(x_\e(t),I_\e(t)\big) = O(\e).$
Hence, combining with the latter equation we obtain \begin{align*}
    &R\big(x_\e(t),I_\e(t)\big) + O(\e) - R\big(\overline{x}(t),\mathcal{I}(x_\e(t)\big) = R\big(x_\e(t), \mathcal I(x_\e(t))\big)-R\big(\overline{x}(t), \mathcal I(x_\e(t))\big)\end{align*}
    Thus, we have 
    \begin{align*}
    R\big(\overline{x}(t),I_\e(t)\big) - R\big(\overline{x}(t),\mathcal{I}(x_\e(t))\big) =& \ R\big(\overline{x}(t),I_\e(t)\big)-R\big(x_\e(t),I_\e(t)\big) +O(\e) \\ &  + R\big(x_\e(t),\mathcal{I}(x_\e(t)\big)-R\big(\overline{x}(t),\mathcal{I}(x_\e(t))\big).
\end{align*}
There exists $ I \in e(I_\e(t),\mathcal{I}(x_\e(t))),$ and $ (c,c') \in [0,1]^2, $ such that  
\begin{align*}
    \frac{\partial}{\partial I} R\big(\overline{x}(t),I)(I_\e(t)-\mathcal{I}(x_\e(t)\big) =& \ O(\e) + \nabla_x R\big(c\overline{x}(t) + (1-c)x_\e(t),I_\e(t)\big)\cdot (\overline{x}(t)-x_\e(t)) \\
     &+ \nabla_x R\big(c'x_\e(t) + (1-c')\overline{x}(t),\mathcal{I}(x_\e(t))\big)\cdot (x_\e(t)-\overline{x}(t)).
\end{align*}
According to the assumption \fer{asrDi}, we have $-\underline{K}_2 \leq \frac{\partial}{\partial I} R\leq - \overline{K}_2 $. Moreover, $t \mapsto \nabla R(\overline{x}(t))$, and $t \mapsto \nabla R(x_\e(t))$ are bounded because $\overline{x}$ and $x_\e$ are in a compact set for $t$ in $[0,T].$\\ 
Finally, we obtain  $$ \vert I_\e(t) - \mathcal I(x_\e(t))\vert \leq C \vert x_\e(t) - \overline x(t)\vert + C\e.$$\\
\textbf{(ii)} To prove the expansions, we will use the following Lemma.
\begin{lemma}
\label{lemma2}
Under assumptions \fer{aspsi}--\fer{asuIni}, for $ \delta \in ]0,T]$ small enough, we have
    $$\|v_\e - u\|_{L^{\infty}_t W^{4,\infty}_x([0,\delta]\times \R^d)} \leq C \| \mathcal{I}(\overline{x}(\cdot))-I_\e\|_{L^\infty([0,\delta])} \delta + C\e \delta.$$ 
 \end{lemma}
\textbf{Proof of Lemma \ref{lemma2}.}
\\\\ \textbf{(a) Let us show that $\| v_\e - u \|_{L^\infty([0,\delta]\times \R^d)}  \leq C \| I_\e -I \|_\infty \delta +C\e\delta$.}\\
To prove this inequality, we define $r := v_\e - u$. According to \fer{system-em} and \fer{system-u}, $r$ satisfies the equation: 
\beq
\label{eqR}
\partial_t r = \e \Delta v_\e + (\nabla u + \nabla v_\e)\cdot \nabla r + R\big(x,I_\e(t)\big)-R\big(x,I(t)\big).
\eeq
We have already proven that $D^2 v_\e$ is uniformly bounded in Theorem \ref{thm:approx}, thus $\e\Delta v_\e = O(\e)$.
Then, we can write:
$$\partial_t r = (\nabla u + \nabla v_\e)\cdot \nabla r + R \big(x,I_\e(t)\big)-R \big(x,I(t)\big) + O(\e).$$
We study this PDE using characteristics. The characteristics satisfy:
$$\dot \gamma (t) = -\nabla u(t,\gamma(t)) - \nabla v_\e(t,\gamma(t)).$$\\
Let $(t_1,x_1)$ be in $ [0,\delta] \times \R^d$. The gradients $\nabla u$ and  $\nabla v_\e $ are Lipschitz with respect to their state variable. Hence, we can apply the Cauchy-Lipchitz Theorem. There exists a unique global characteristic $\gamma(t)$ for $t\in [0,t_1]$ such that $\gamma(t_1)=x_1$. 
Multiplying $\dot \gamma(t)$ by $\nabla r (t,\gamma(t))$, we obtain 
$$ \nabla r \big(t,\gamma(t)\big)\cdot \dot \gamma(t) = \Big(-\nabla u\big(t,\gamma(t)\big) - \nabla v_\e\big(t,\gamma(t)\big)\Big)\cdot \nabla r \big(t,\gamma(t)\big).$$
Thus with \fer{eqR}, we obtain
$$ \frac{\mathrm{d}}{\mathrm{d}t}r\big(t,\gamma(t)\big) = O(\e) + R\big(x,I_\e(t)\big)-R\big(x,I(t)\big).$$
We integrate and use the initial condition $r\big(0,\gamma(0)\big) =O(\e) $, to obtain $$ r\big(t_1,\gamma(t_1)\big) = r(t_1,x_1) = \int_0^{t_1} R\big(x,I_\e(s)\big)-R\big(x,I(s)\big)\mathrm{d}s + O(\e)t_1.$$
Finally, we find $$ \forall (t_1,x_1)\in [0,\delta]\times \R^d, \quad \vert r(t_1,x_1)\vert \leq C \| I_\e -I \|_{L^\infty([0,\delta])} \delta +C\e\delta.$$
This is the inequality we wanted $\| v_\e - u \|_{L^\infty([0,\delta]\times \R^d)}  \leq C \| I_\e -I \|_{L^\infty([0,\delta])} \delta +C\e\delta$.\\\\
\textbf{(b) Let us prove that $\| \nabla v_\e -\nabla u \|_{L^\infty([0,\delta]\times \R^d)}  \leq C \| I_\e -I \|_{L^\infty([0,\delta])} \delta +C\e\delta$}.
\\ We differentiate with respect to $e_i$ the equation \fer{eqR}, and we obtain
$$ \partial_t \partial_i r = \e \Delta \partial_i v_\e + (\nabla\partial_i v_\e+\nabla\partial_i u)\cdot \nabla r + (\nabla u + \nabla v_\e)\cdot \nabla \partial_i r + \partial_i R(x,I_\e) - \partial_i R(x,I).$$
Multiplying by $\partial_i r$ and summing over all $i$, we have 
\begin{align*}
    \sum_{i=1}^d \partial_t \partial_i r \cdot \partial_i r =& \ \e \nabla \Delta v_\e \cdot \nabla r +\sum_{i=1}^d (\nabla\partial_i v_\e+\nabla\partial_i u)\cdot \nabla r \partial_i r+ \sum_{i=1}^d (\nabla v_\e + \nabla u)\cdot \nabla \partial_i r \partial_i r \\
    &+ \sum_{i=1}^d (\partial_i R(x,I_\e) - \partial_i R(x,I))\partial_i r.
\end{align*}
Thus, we deduce
\begin{align*}
    \sum_{i=1}^d \partial_t \partial_i r \cdot \partial_i r \leq & \ \e\vert \nabla \Delta v_\e \vert \vert  \nabla r \vert + \sum_{i=1}^d \vert\nabla\partial_i v_\e+\nabla\partial_i u \vert \vert \partial_i r \vert \vert \nabla r \vert + \vert \nabla R(x,I_\e)-\nabla R(x,I) \vert \vert \nabla r \vert\\
    &+ \sum_{i=1}^d (\nabla v_\e + \nabla u)\cdot \nabla \partial_i r \partial_i r.
\end{align*}
We know that $D^2 v_\e $, $D^2 u$, $D^3 v_\e$ are uniformly bounded thanks to \fer{eq.D2u}, and \fer{nD3u}, then we divide by $\vert \nabla r \vert$,\\ and obtain 
$$\partial_t \vert \nabla r \vert \leq C\e + C \vert \nabla r \vert + (\nabla v_\e + \nabla u)\cdot\nabla \vert \nabla r \vert + \vert \nabla R(x,I_\e)-\nabla R(x,I) \vert .$$
Finally, we have \beq \label{ineqR}\partial_t \vert \nabla r \vert \leq C\e + C \vert \nabla r \vert + (\nabla v_\e + \nabla u)\cdot\nabla \vert \nabla r \vert + C \vert I_\e(t)-I(t) \vert .\eeq
As in the last point, we study the characteristics of the equation.
The characteristics satisfy: $$\dot \gamma (t) = -\nabla u(t,\gamma(t)) - \nabla v_\e(t,\gamma(t)).$$ 
Let $(t_1,x_1)$ be in $ [0,\delta] \times \R^d$, according to the Cauchy-Lipshitz Theorem, there exists a unique characteristic such as $\gamma(t_1)=x_1.$ Moreover, $\gamma$ is well-defined in $[0,t_1].$
By definition of $\dot \gamma$, we have
$$\nabla \vert \nabla r \vert (t,\gamma(t))\cdot \dot \gamma(t) =  -(\nabla u(t,\gamma(t)) +\nabla v_\e(t,\gamma(t)))\cdot \nabla \vert \nabla r \vert (t,\gamma(t)).$$
We add the equation \fer{ineqR} to the latter equation, and we obtain:  $$ \frac{\mathrm{d}}{\mathrm{d}t} \vert \nabla r \vert (t,\gamma(t)) \leq C \vert \nabla r \vert (t,\gamma(t))+ C \| I_\e-I\|_{L^\infty([0,\delta])} + C\e.$$
Using the Grönwall's Lemma, we obtain the inequality.
\\\\ \textbf{(c) We now prove $\|D^2 v_\e -D^2 u \|_{L^\infty([0,\delta]\times \R^d)}  \leq C \| I_\e -I \|_{L^\infty([0,\delta])} \delta +C\e\delta$.}\\
We differentiate two times with respect to $\xi$ the equation \fer{eqR}, we find
\begin{align*}
\partial_t \partial_{\xi\xi}^2 r =& \ \e \Delta \partial_{\xi\xi}^2 v_\e + (\nabla \partial_{\xi\xi}^2v_{\e} + \nabla \partial_{\xi\xi}^2 u\cdot \nabla r +  2(\nabla \partial_\xi v_{\e} + \nabla \partial_\xi u)\cdot \nabla \partial_\xi r + \partial_{\xi\xi}^2 R(x,I_\e)- \partial_{\xi\xi}^2 R(x,I)\\ & + (\nabla v_\e + \nabla u)\cdot \nabla \partial^2_{\xi\xi} r.\end{align*} 
Thanks to the estimation \eqref{nD3u} in Theorem \ref{t2.2}, we know that $D^4 v_\e$ is uniformly bounded. Thus, $\e \Delta \partial_{\xi\xi}^2 v_\e =O(\e)$.
With similar arguments, we obtain
\beq \label{r2}\partial_t \vert \partial_{\xi\xi}^2 r \vert \leq C \|I_\e-I\|_{L^\infty([0,\delta])} \delta + C \|\partial_{\xi\xi}^2 r \|_{L^\infty([0,\delta]\times \R^d)} + (\nabla v_\e+\nabla u)\cdot \nabla \vert \partial_{\xi\xi}^2 r \vert + C\| I-I_\e\vert \|_{L^\infty([0,\delta])}+ C\e. \eeq
The characteristics follow  $$\dot \gamma (t) = -\nabla u(t,\gamma(t)) - \nabla v_\e(t,\gamma(t)).$$ 
Let $(t_1,x_1)$ be in $[0,\delta] \times \R^d$, and according to the Cauchy-Lipshitz Theorem, there exists a unique trajectory $\gamma$ such that $\gamma(t_1)=x_1$, and it is well defined in $[0,t_1].$
By definition of $\dot \gamma$, we have $$ \nabla \vert\partial_{\xi\xi}^2 r(t,\gamma(t))\vert \cdot \gamma(t) = -(\nabla v_\e(t,\gamma(t))+\nabla u(t,\gamma(t)))\cdot \nabla \vert \partial_{\xi\xi}^2 r(t,\gamma(t))\vert .$$\\
We use the equation \fer{r2}, to obtain
$$\frac{\mathrm{d}}{\mathrm{d}t} \vert \partial_{\xi\xi}^2 r(t,\gamma(t)) \vert \leq C \| I-I_\e \|_{L^\infty([0,\delta])}+  C \|\partial_{\xi\xi}^2 r\|_{L^\infty([0,\delta]\times \R^d)} + C\e.$$
We integrate the latter line in $[0,\delta]$ and deduce that $$\vert \partial_{\xi\xi}^2 r(t_1,x_1) \vert \leq C\delta \| I-I_\e \|_{L^\infty([0,\delta])} + C \|\partial_{\xi\xi}^2 r\|_{L^\infty([0,\delta]\times \R^d)}\delta + C\e \delta .$$
Hence $\| \partial_{\xi\xi}^2 r \| _{L^{\infty}([0,\delta]\times \R^d)}\leq C\delta \| I-I_\e \|_{L^\infty([0,\delta])}+ C \|\partial_{\xi\xi}^2 r\|_{L^{\infty}([0,\delta]\times \R^d)} \delta + C\e \delta,$
and for $\delta $ sufficiently small, we obtain the wanted result. \\\\
\textbf{(d)} We prove $\|D^k v_\e -D^k u \|_{L^\infty([0,\delta]\times \R^d)}  \leq C \| I_\e -I \|_{L^\infty([0,\delta])} \delta +C\e\delta$, for $k=3,4,$ with similar arguments.
\begin{flushright}
    $\square$
\end{flushright}
\textbf{(iii)} We now show that  $\| x_\e-\overline{x} \|_{L^\infty([0,\delta])}\leq C \delta\| x_\e-\overline{x} \|_{L^\infty([0,\delta])} + C\e\delta +C\e$. \\\\
According to the equations \fer{system-em} and \fer{system-u}, $x_\e$ and $\overline{x}$ satisfy the following ODE:
\begin{align*}
\dot{ x}_\e(t) &= \left( -D^2 v_\e \big(t,  x_\e(t)\big) \right)^{-1}  \left[ \nabla_x R\big( x_\e(t),I_\e(t) \big) +\e D^3 v_\e \right], \\
    \dot{\overline{ x}}(t) &= \left( -D^2u \big(t,  \overline{x}(t)\big) \right)^{-1} \nabla_x R\big( \overline{x}(t),I(t) \big).
\end{align*}
Moreover $x_\e(0)=\overline{x}(0) +O(\e) = x_0 + O(\e)$, thus 
\begin{align*}
    \overline{x}(t)-x_\e(t) = & \ \int_0^t \Big[\big( -D^2u \big(s,  \overline{x}(s)\big) \big)^{-1}   \nabla_x R\big( \overline{x}(s),I(s) \big) \\
   & -\left( -D^2 v_\e \big(s,  x_\e(s)\big) \right)^{-1}  \left[ \nabla_x R\big( x_\e(s),I_\e(s) \big) +\e \nabla \Delta v_\e \right] \Big] \mathrm{d}s +O(\e), \end{align*}
   and hence, 
   \begin{align*}
   \vert\overline{x}(t)-x_\e(t) \vert  \leq &\ \e \int_0^t \vert \big(D^2v_\e(s,x_\e(s))\big)^{-1} \nabla \Delta v_\e(s,x_\e(s)) \vert \mathrm{d}s\\
   & + \int_0^t \big \vert \big[\big(-D^2u(s,\overline{x}(s)\big)^{-1} - \big(-D^2 v_\e(s,x_\e(s))^{-1}\big)\big]\nabla_x R\big(\overline{x}(s), I(s)\big)\big\vert \mathrm{d}s \\
   & + \int_0^t\big \vert \big(-D^2 v_\e(s,x_\e(s))\big)^{-1} [\nabla_x R\big(x_\e(s),I_\e(s)\big) - \nabla_x R \big(\overline{x}(s),I(s)\big)]\big\vert \mathrm{d}s +O(\e).
\end{align*}
We study the three integrals of the right-hand side independently, we define:
\begin{align*}
    A_1 &:=\e \int_0^t \big\vert \big(D^2v_\e(s,x_\e(s))\big)^{-1} \nabla \Delta v_\e(s,x_\e(s)) \big \vert \mathrm{d}s, \\
    A_2 &:=  \int_0^t \big \vert \big[\big(-D^2u(s,\overline{x}(s)\big)^{-1} - \big(-D^2 v_\e(s,x_\e(s))^{-1}\big)\big]\nabla_x R\big(\overline{x}(s), I(s)\big)\big\vert \mathrm{d}s,\\
    A_3 &:= \int_0^t\big \vert \big(-D^2 v_\e(s,x_\e(s))\big)^{-1} [\nabla_x R\big(x_\e(s),I_\e(s)\big) - \nabla_x R \big(\overline{x}(s),I(s)\big)]\big\vert \mathrm{d}s.
\end{align*}
-- According to the results \fer{eq.D2u} and \fer{nD3u} on $D^2 v_\e$ and $D^3 v_\e$, we have
$$A_1 \leq \frac{\e \delta C(T)}{\overline{L}_2} .$$
-- We now estimate $A_2$.\\
We can write the following equality
$$ \big(-D^2u(s,\overline{x}(s))\big)^{-1} - (-D^2 v_\e(s,x_\e(s)))^{-1} = $$ $$\big(-D^2u(s,\overline{x}(s))\big)^{-1} \big[D^2u(s,\overline{x}(s))-D^2 v_\e(s,x_\e(s))\big] \big(-D^2 v_\e(s,x_\e(s))\big)^{-1},$$
and we use the result of Lemma \ref{lemma2}: 
$$ \| v_\e - u \|_{L^\infty([0,\delta])\times W^{2,\infty}( \R^d)}\leq C \delta \|I_\e - \mathcal I \big(\overline{x}(\cdot)\big) \|_{L^\infty([0,\delta])} + C\e \delta.$$
We write the following triangular inequality: $$\vert I_\e(t) - \mathcal I (\overline{x}(t)) \vert \leq \vert I_\e(t) -\mathcal I (x_\e(t)) \vert + \vert \mathcal I (x_\e(t)) - \mathcal I (\overline{x}(t)) \vert .$$
We already know from (i) that  $$ \vert I_\e(t) -\mathcal I (x_\e(t)) \vert \leq C\e + C \vert x_\e(t)-\overline{x}(t) \vert, $$
and we could show similarly, using $R\big(x_\e(t),\mathcal I (x_\e(t))\big) = 0 = R\big(\overline{x}(t), \mathcal{I}(\overline{x}(t)\big)$, that $\vert \mathcal I \big(x_\e(t)\big) - \mathcal I \big(\overline{x}(t)\big) \vert \leq C \vert x_\e(t)-\overline{x}(t) \vert$. We obtain \beq \label{ineq1} |I_\e(t)-\mathcal I(\bar x(t))|\leq C\e+C|x_\e(t)-\bar x(t)|.\eeq
Thus for $\delta <1$, we have: $$A_2 \leq \delta \| \nabla R(\overline{x}(\cdot),I(\cdot))) \|_{L^\infty ([0,T])} \times (C\delta  \| x_\e - \overline{x} \| _{L^\infty([0,\delta])}+ C\e\delta)
\leq C\e \delta + C\delta \| x_\e - \overline{x}\| _{L^\infty([0,\delta])}.
$$
-- We estimate $A_3$.\\
Using assumptions \fer{asrD2}, \fer{asr23} and \fer{ineq1}, we obtain
\begin{align*}
    \vert \nabla_x R (x_\e(s),I_\e(s)) -\nabla_x R (\overline{x}(s),I(s)) \vert \leq & \ \vert \nabla_x R (x_\e(s),I_\e(s) - \nabla_x R (\overline{x}(s),I_\e(s)) \vert \\
    & + \vert \nabla_x R (\overline{x}(s),I_\e(s))  - \nabla_x R (\overline{x}(s),I(s)) \vert \\
     \leq & \ \sup_{t\in [0,T]} \| D^2 R(\cdot,s) \|_\infty \| x_\e - \overline{x} \| _{L^\infty([0,\delta])} \\ &+ C \sup_{x,i} |\frac{\partial^2}{\partial I \partial x_i} R| |I_\e(s) - I(s)| \\
     \leq & \  C\e + C  \| x_\e - \overline{x}\| _{L^\infty([0,\delta])}.
\end{align*}   
Thus, we proved that $ \| x_\e - \overline{x}\| _{L^\infty([0,\delta])} \leq C\e \delta + C\delta \| x_\e - \overline{x}\| _{L^\infty([0,\delta])} + C\e.$
Then $(1-C \delta)\| x_\e - \overline{x}\| _{L^\infty([0,\delta])} \leq C \e \delta+C\e,$
and for $\delta < \frac{1}{2C}$, we have an estimate on $x_\e$ for $t$ in $[0,\delta]$ $$\| x_\e - \overline{x}\| _{L^\infty([0,\delta])} \leq C \e .$$ We can iterate the argument $\lceil \frac{T}{\delta} \rceil$ times to obtain the result $$ \| x_\e - \overline{x}\| _{L^\infty([0,T])} \leq C \e.$$
We combine this inequality with Lemma \ref{lemme1}, and Lemma \ref{lemma2}, to conclude the proof of Proposition \ref{item1}
$$
\|I_\e  -I\|_{L^\infty([0,T])}\leq C\e,\ \| x_\e -\overline x \|_{L^\infty([0,T])} \leq C\e, \
\| u_\e - u- \e \log (\f{r}{ \e^{\f d 2}}) \|_{L^\infty_tW^{4,\infty}_x([0,T]\times \R^d)}\leq C \e. $$
\begin{flushright}
    $\square$
\end{flushright}
\section{\texorpdfstring{The full expansion of $I_\e$ and $u_\e$}{}}
In this section, we prove Theorem \ref{thm:approx}. We first need an approximation of order $1$ for $I_\e$. To this end, we use the result of Proposition \ref{item1} and we use similar ideas as in Section \ref{s3} to prove the following Lemma \ref{lemmeK}.
\begin{lemma}
\label{lemmeK}
Under assumptions \fer{aspsi}--\fer{asuIni}, there exists a constant $C(T) = C>0$, independant of $\e$, such that we have
$$
\| I_\e - \mathcal I\big(x_\e(\cdot)\big) -\e K\|_{L^\infty([0,T])}\leq C\e^2,
$$
with
$$
\begin{array}{rl}
 K(t)=& -
 \f{1 }{\f{\p }{\p I}R(\overline x(t),I(t))}
 \left[\left( \f{\nabla \psi(\overline x(t))}{\psi(\overline x(t))} \nabla R\big( \overline x(t), I(t)\big )+  \f {1}{2} D^2 R \big(\overline x(t), I(t) \big) \right) \left| D^2 u(t,\overline x(t) ) \right|^{-1}\right.
 \\[2mm]
& +   \f {1}{2} D^3 u(t,\overline x(t)) \nabla R \big ( \overline x(t), I(t) \big )  \left| D^2 u(t,\overline x(t) ) \right|^{-2}
\left.- \f{\nabla_x R(\overline x(t),I(t)) \left(-D^2 u(t,\overline x(t)) \right)^{-1}\nabla R \big(\overline x(t),I(t)\big)}{I(t) \f{\p }{\p I}R(\overline x(t),I(t))}\right].
\end{array}
$$
\end{lemma}
\textbf{Proof.}
We give the proof of Lemma \ref{lemmeK} in the appendix \ref{A_K}.\\\\
 To obtain an asymptotic expansion of $x_\e$ and $v_\e$, we study the following functions: $y_\e := \frac{x_\e-\overline x}{\e}$ and $w_\e := \frac{v_\e-u}{\e}.$
\subsection{\texorpdfstring{Equations satisfied by $y_\e$ and $w_\e$}{}}
We start by showing that $y_\e$ and $w_\e$ satisfy the system of Lemma \ref{lemme6}.
\begin{lemma}
\label{lemme6}
    In $[0,T]$, $(w_\e,y_\e)$ satisfy the system
    \begin{System}
    \label{systeme_e}
        \dot{y}_\e(t) = \Tilde{f}(t) + \Tilde{g}(t) \nabla \mathcal{I}(\overline{x}(t))\cdot y_\e(t) 
          +[-D^2 u(t,\overline{x}(t))]^{-1} D^2 w_\e(t,\overline{x}(t)) [-D^2 u(t,\overline{x}(t))]^{-1} \nabla R(\overline{x}(t),I(t)) \\
        \qquad \quad \ \ +[-D^2 u(t,\overline{x}(t))]^{-1} D^3 u(t,\overline{x}(t))y_\e(t) [-D^2 u(t,\overline{x}(t))]^{-1} \nabla R(\overline{x}(t),I(t))+ \Tilde{h}(t)y_\e(t) + O(\e) \\
         \partial_t w_\e = 2\nabla u \cdot \nabla w_\e + \frac{\partial}{\partial I} R(x,I(t)) \nabla \mathcal{I}(\overline{x}(t))\cdot y_\e(t) + \Delta u + \frac{\partial}{\partial I} R(x,I(t))K(t)+O(\e)
    \end{System}
     with initial conditions: $y_\e(0) = 0,\ w_\e(0,x) = 0$,
     and \begin{align*}
         \Tilde{f}(t) &:= \nabla \Delta u(t,\overline{x}(t)) + [-D^2 u(t,\overline{x}(t)]^{-1} \frac{\partial}{\partial I} \nabla R(\overline{x}(t),I(t))K(t),\\ 
         \Tilde{g}(t) &:= \frac{\partial}{\partial I} \nabla R(\overline{x}(t),I(t)),\
         \\ \Tilde{h} (t) &:= [-D^2 u(t,\overline{x}(t))]^{-1} D^2 R(\overline{x}(t),I(t)). 
     \end{align*} Moreover, we have estimates for the derivatives of $w_\e:$
     $$\partial_t \partial_i w_\e = \Delta \partial_i u + 2 \nabla \partial_i u \cdot \nabla w_\e + 2 \nabla  u \cdot \nabla \partial_i w_\e + \partial_I \partial_i R(x,I(t))\Big(\mathcal{I}(\overline x(t)) \cdot y_\e(t) +K(t)\Big)+O(\e),$$
and 
\begin{align*}
    \partial_t \partial_{i,j}^2 w_\e =&  \ \Delta \partial_{i,j}^2 u + 2 \nabla \partial_{i,j}^2 u \cdot \nabla w_\e + 2 \nabla \partial_j u \cdot \nabla \partial_i w_\e + 2 \nabla u\cdot \nabla \partial_{i,j}^2 w_\e \\ &+  \partial_I \partial_{i,j}^2 R(x,I(t))\Big(\mathcal{I}(\overline x(t)) \cdot y_\e(t)+K(t)\Big) +O(\e) \nonumber .\end{align*}
\end{lemma}
\textbf{Proof.} Using the equations \fer{system-em} and \fer{system-u}, we have
\begin{align}
\label{x_e-x}
    \dot{x}_\e(t)-\dot{\overline{ x}}(t) =& \ \e \nabla \Delta v_\e(s,x_\e(s)) + ([-D^2v_\e(s,x_\e(s))]^{-1} - [-D^2 u(s,\overline x (s))]^{-1} )\nabla R(x_\e(s),I_\e(s))\\
    &- [-D^2 u(s,\overline x (s))]^{-1} [\nabla R(\overline{x}(t),I(s))- \nabla R(x_\e(s),I_\e(s))].\nonumber
\end{align}
We know from Lemma \ref{lemma2} that $ \e  \nabla \Delta v_\e(s,x_\e(s))  = \e  \nabla \Delta u(s,\overline{x}(s)) +O(\e^2)$.
\\\\
We estimate the last term of the right-hand side of \fer{x_e-x} 
$$ \nabla R(\overline{x}(t),I(s))- \nabla R(x_\e(s),I_\e(s)) = \nabla R(\overline{x}(t),I(s))- \nabla R (\overline{x}(s),I_\e(s)) + \nabla  R (\overline{x}(s),I_\e(s))- \nabla R(x_\e(s),I_\e(s)),$$
and use a Taylor expansion to obtain
\begin{align*}
    - \nabla R(\overline{x}(s),I(s)) + \nabla R (\overline{x}(s),I_\e(s))  =& \frac{\partial}{\partial I} \nabla R(\overline{x}(t),I(s))(I_\e(s)-I(s)) +  \frac{\partial^2}{\partial I} \nabla R(\overline{x}(t),I(s))(I_\e(s)-I(s))^2 \\
    & + o((I_\e(s)-I(s))^2).
\end{align*}
According to Proposition \ref{item1}, we already know that $I_\e(t) = I(t) +O(\e)$, thus $(I_\e(s)-I(s))^2 = O(\e^2)$,\\\\
and we find $-\nabla R(\overline{x}(s),I(s)) + \nabla R (\overline{x}(s),I_\e(s)) = \frac{\partial}{\partial I} \nabla R(\overline{x}(t),I(s))(I_\e(s)-I(s))+O(\e^2)$.
\\\\We also know from Proposition \ref{item1} that $x_\e = \overline{x} +O(\e)$,  which gives us\\\\
 $-\nabla R(\overline{x}(t),I_\e(t))+\nabla R(x_\e(t),I_\e(t)) = 
D^2 R(\overline{x}(t),I_\e(t))(\overline{x}(t)-x_\e(t)) +O(\e^2).$
\\\\ 
We obtain the following equation
\begin{align*}
    \dot{x}_\e(t)-\dot{\overline{ x}}(t) =& \  \e  \nabla \Delta u(s,\overline{x}(s)) + ([-D^2v_\e(s,x_\e(s))]^{-1} - [-D^2 u(s,\overline x (s))]^{-1} )\nabla R(x_\e(s),I_\e(s)) \\
    & + [-D^2u(s,\overline{x}(s))]^{-1} \frac{\partial}{\partial I} \nabla R (\overline{x}(s),I(s))(I_\e(s)-I(s)) \nonumber \\ 
    &+  [-D^2 u (s,\overline{x}(s))]^{-1} D^2 R(\overline{x}(s),I_\e(s))(\overline{x}(s)-x_\e(s)) +O(\e^2). \nonumber
\end{align*}
We proved in Lemma \ref{lemmeK} that $$ I_\e(t) = \mathcal{I}(x_\e(t)) +K(t)\e +O(\e^2).$$
Then, we write
\begin{align*}
I_\e(t)-I(t) &= I_\e(t) - \mathcal{I}(x_\e(t)) +  \mathcal{I}(x_\e(t)) - I(t)\\
&= \e K(t) + O(\e^2) + \mathcal{I} (x_\e(t))- I(t)\\
&= \e K(t) +  \nabla \mathcal I (\overline{x}(t)) \cdot (x_\e(t)-\overline{x}(t))  + \f 12 D^2\mathcal I (\overline{x}(t)) \cdot (x_\e(t)-\overline{x}(t))^2 + O(\e^2).
\end{align*}
We obtain a first expansion of $I_\e(t)$ $$  I_\e(t) = I(t) + \e K(t) + \nabla \mathcal I (\overline{x}(t)) \cdot (x_\e(t)-\overline{x}(t))+ O(\e^2).$$
We deduce that \begin{align*}
    [-D^2u(s,\overline{x}(s))]^{-1} \frac{\partial}{\partial I} \nabla R (\overline{x}(s),I(s))(I_\e(s)-I(s)) = \e  [-D^2 u (s,\overline{x}(s)]^{-1} \frac{\partial}{\partial I} \nabla R(\overline{x}(s),I(s)) K(s) \\
    + [-D^2 u (s,\overline{x}(s))]^{-1} \frac{\partial}{\partial I} \nabla R(\overline{x}(s),I(s)) \nabla \mathcal{I}(\overline{x}(s))\cdot (x_\e(s)-\overline{x}(s))+O(\e^2).
\end{align*} 
We next write a Taylor expansion for $D^2 R(\overline{x}(s),I_\e(s))$,
\begin{align*}
    D^2 R (\overline{x}(s),I_\e(s)) & = D^2 R(\overline{x}(s),I(s)) + \e \frac{\partial}{\partial I} D^2 R(\overline{x}(s),I(s))(K(s)+\nabla \mathcal{I}(\overline{x}(s)) \cdot y_\e(s))+O(\e^2)\\
    &= D^2 R(\overline{x}(s),I(s))+ \e \frac{\partial}{\partial I} D^2 R(\overline{x}(s),I(s))K(s)+ O(\e^2).
\end{align*}
We now only need to write a Taylor expansion for the second term of the right-hand side of \fer{x_e-x} 
\begin{align*}
&[-D^2 v_\e(t,x_\e(t))]^{-1}- [-D^2 u (t,\overline{x}(t))]^{-1} \\= &  \ [-D^2 v_\e(t,x_\e(t))]^{-1} [D^2 v_\e(t,x_\e(t)) -D^2 u(t,\overline{x}(t))][-D^2 u(t,\overline{x}(t))]^{-1}. \end{align*}
We estimate the following expression
\begin{align*}
    D^2 v_\e(t,x_\e(t)) -D^2 v_\e(t,\overline{x}(t)) + D^2 v_\e(t,\overline{x}(t)) - D^2 u(t,\overline{x}(t)) =& \ D^3 v_\e(t,\overline{x}(t))(\overline{x}(t)-x_\e(t)) + O(\e^2)\\& + D^2(v_\e-u)(t,\overline{x}(t)) \\
    =& \ D^2(v_\e-u)(t,\overline{x}(t)) \\&+\e D^3 v_\e(t,\overline{x}(t))y_\e(t) +O(\e^2).
\end{align*}
Finally, $y_\e$ satisfies the ODE \\
\begin{align*}
           \dot{y}_\e(t) =& \ \Tilde{f}(t) + \Tilde{g}(t)\nabla \mathcal{I}\cdot y_\e(t) + \Tilde{h}(t)y_\e(t) +  [-D^2 u(t,\overline{x}(t))]^{-1} D^2 w_\e(t,\overline{x}(t)) [-D^2 u(t,\overline{x}(t))]^{-1} \nabla R(\overline{x}(t),I(t))\\
         & + [-D^2 u(t,\overline{x}(t))]^{-1} D^3 u(t,\overline{x}(t))y_\e(t) [-D^2 u(t,\overline{x}(t))]^{-1} \nabla R(\overline{x}(t),I(t))+O(\e).         
\end{align*}
We now look for an equation on $w_\e,$ to this end we use systems \fer{system-em}, \fer{system-u}. By subtraction, we obtain $ \partial_t (v_\e-u) = \e \Delta v_\e + (\nabla v_\e+\nabla u)\cdot \nabla (v_\e-u) + R(x,I_\e(t))-R(x,I(t)).$\\\\
We previously proved in Lemma \ref{lemma2}, that for $\delta$ small enough, we have $\Delta v_\e = \Delta u + O(\e)$ in $[0,T]\times \R^d$, hence $\e \Delta v_\e = \e \Delta u + O(\e^2)$ in $[0,T]\times \R^d$.
\\\\ We use again the expression: $I_\e(t) -I(t) = \e K(t) + \nabla \mathcal{I}(\overline{x}(t))\cdot(x_\e(t) -\overline{x}(t))+O(\e^2)$, and we obtain 
\begin{align*}
    R(x,I_\e(t)) -R(x,I(t))  =& \ \frac{\partial}{\partial I} R(x,I(t))\cdot [\e K(t) + \nabla \mathcal{I}(\overline{x}(t))\cdot (x_\e(t)-\overline{x}(t)) +O(\e^2)]\\
    &+ \f 12 \frac{\partial^2}{\partial I^2} R(x,I(t)) (\overline{x}(t)-x_\e(t))^2 + O(\e^2).
\end{align*}
We deduce that 
$$R(x,I_\e(t)) -R(x,I(t)) = \e \frac{\partial}{\partial I} R(x,I(t))\cdot  K(t) + \frac{\partial}{\partial I} R(x,I(t))\cdot \nabla \mathcal{I}(\overline{x}(t))\cdot (x_\e(t)-\overline{x}(t)) + O(\e^2).$$
We obtain an equation for $w_\e$: 
$$\partial_t w_\e = \Delta u + \frac{\partial}{\partial I} R(x,I(t))K(t) + (\nabla v_\e + \nabla u) \cdot \nabla w_\e + \frac{\partial}{\partial I} R(x,I(t)) \nabla \mathcal{I}(\overline{x}(t))\cdot y_\e(t) + O(\e).$$
We also know from Lemma \ref{lemma2}, that in $[0,T]\times\R^d$, $\nabla v_\e = \nabla u + O(\e).$
Thus, $\nabla v_\e \cdot \nabla (v_\e-u) = \nabla u \cdot \nabla (v_\e-u) + O(\e^2)$, and  $\nabla v_\e \cdot \nabla w_\e = \nabla u \cdot \nabla w_\e + O(\e). $\\\\
Finally, $w_\e$ satisfies the following transport equation 
$$ \partial_t w_\e = 2\nabla u \cdot \nabla w_\e + \frac{\partial}{\partial I} R(x,I(t)) \nabla \mathcal{I}(\overline{x}(t))\cdot y_\e(t) + \Delta u + \frac{\partial}{\partial I} R(x,I(t))K(t)+ O(\e).$$ 
With similar arguments, we obtain equations for the derivatives of $w_\e$. We have
$$\partial_t \partial_i w_\e = \Delta \partial_i u + 2 \nabla \partial_i u \cdot \nabla w_\e + 2 \nabla  u \cdot \nabla \partial_i w_\e + \partial_I \partial_i R(x,I(t))\Big(\mathcal{I}(\overline x(t)) \cdot y_\e(t) +K(t)\Big)+O(\e),$$
and 
\begin{align*}
    \partial_t \partial_{i,j}^2 w_\e =&  \ \Delta \partial_{i,j}^2 u + 2 \nabla \partial_{i,j}^2 u \cdot \nabla w_\e + 2 \nabla \partial_j u \cdot \nabla \partial_i w_\e + 2 \nabla u\cdot \nabla \partial_{i,j}^2 w_\e \\ &+  \partial_I \partial_{i,j}^2 R(x,I(t))\Big(\mathcal{I}(\overline x(t)) \cdot y_\e(t)+K(t)\Big) +O(\e).\end{align*}
The couple $(y_\e,w_\e)$ satisfies the system written in Lemma \ref{lemme6} and this ends the proof.
\begin{flushright}
    $\square$
\end{flushright}
To finish the proof of Theorem \ref{thm:approx} (ii), we need to prove that $y_\e$ and $w_\e$ converge to a limit in $L_t^\infty W^{2,\infty}_x([0,T]\times \R^d)$, when $\e$ approaches $0$. To this end, we prove the Lemma \ref{dernierlemme}.
\begin{lemma}
\label{dernierlemme}
    The couple $(y_\e,w_\e)_\e$ uniformly converges in $ L_t^\infty W^{2,\infty}_x([0,T]\times \R^d)$ to the unique solution $(y,w)$ of the following system:
 \begin{System}
 \label{syst_limit}
        \dot{y}(t) = \Tilde{f}(t) + \Tilde{g}(t)\cdot y(t) + \Tilde{h}(t)y(t) +[-D^2 u(t,\overline{x}(t))]^{-1} D^2 w(t,\overline{x}(t)) [-D^2 u(t,\overline{x}(t))]^{-1} \nabla R(\overline{x}(t),I(t)) \\ \quad \quad \ \ \
        + [-D^2 u(t,\overline{x}(t))]^{-1} D^3 u(t,\overline{x}(t))y(t) [-D^2 u(t,\overline{x}(t))]^{-1} \nabla R(\overline{x}(t),I(t)) \\
         \partial_t w = 2\nabla u \cdot \nabla w + \frac{\partial}{\partial I} R(x,I(t)) \nabla \mathcal{I}(\overline{x}(t))\cdot y(t) + \Delta u + \frac{\partial}{\partial I} R(x,I(t))K(t).
\end{System}
with initial conditions: $y(0) = 0,\  w(0,x) = 0.$
\end{lemma}
\textbf{Proof.} We first prove that $(y_\e,w_\e)$ converges in $L_t^\infty W^{2,\infty}_x([0,T]\times \R^d)$ to $(y,w)$ solution of \fer{syst_limit}, and next we will prove the uniqueness of \fer{syst_limit}.\\\\
\textbf{(i)} To this end, we denote $\tilde{y}_\e:= y_\e - y$ and $\tilde{w}_\e := w_\e -w.$ We use \fer{systeme_e} and \fer{syst_limit} to write 
$$\begin{cases}
    \dot {\tilde{y}}_\e(t) &= A(t) \tilde{y}_\e(t) + B(t) D^2 \tilde{w}_\e(t,\overline x(t)) C(t)+E(t)\tilde{y}_\e(t)F(t) +O(\e)\\
    \partial_t \tilde{w}_\e &= 2 \nabla u \cdot \nabla \tilde{w}_\e +G(t,x)\tilde{y}_\e(t)+O(\e)
,\end{cases}$$ with initial conditions: $\tilde{y}_\e(0) =0,\  \tilde{w}_\e(0,x) = 0 $. The maps $A, \ B, \ C, \ D, \ E, \ F, \ G$ are determined in function of the parameter of the problem. We determine $\tilde{w}_\e$ along the characteristics that satisfy:
$$\dot{\gamma}(t) = -2 \nabla u (t,\gamma(t)).$$
According to the Cauchy-Lipshitz Theorem, for $t_1 \in [0,T]$, and for $x_1\in \R^d$, there exists a unique trajectory $\gamma$ defined in $[0,t_1]$ such that $\gamma(t_1)=x_1$.
As before, we obtain
$$ \forall t \in [0,t_1], \quad\frac{\mathrm{d}}{\mathrm{d}t}\vert \tilde{w}_\e\vert (t,\gamma(t)) \leq C \vert \tilde{y}_\e\vert(t)+C\e.$$
We integrate, to find
$$ \forall t\in [0,T], \quad \forall x \in \R^d, \quad |\tilde{w}_\e(t,x)|\leq C\int_0^t |\tilde{y}_\e|(s)\mathrm{d}s+C\e. $$
We take the supremum in space of the last inequality, and obtain
$$ \|\tilde{w}_\e(t,\cdot)\|_{L^\infty(\R^d)} \leq C\int_0^t |\tilde{y}_\e|(s)\mathrm{d}s+C\e.$$
We differentiate the second equation with respect to $x$, and we use that $\nabla G$ and $D^2 u$ are bounded to obtain the inequality
$$ \frac{\mathrm{d}}{\mathrm{d}t} \vert \nabla \tilde{w}_\e \vert (t,\gamma(t)) \leq C \vert \nabla \tilde{w}_\e \vert (t,\gamma(t)) + C\vert \tilde{y}_\e(t) \vert +C\e,$$
then we integrate and  we use that $D^2G$ and $D^3 u$ are bounded to find
$$ \|\nabla\tilde{w}_\e(t,\cdot)\|_{L^\infty(\R^d)}\leq C\int_0^t |\tilde{y}_\e|(s)\mathrm{d}s + C\int_0^t \|\nabla \tilde{w}_\e(s,\cdot)\|_{L^\infty(\R^d)} \mathrm{d}s+C\e.$$
Similarly, we obtain an estimate of the second derivative
$$\frac{\mathrm{d}}{\mathrm{d}t} \vert D^2 \tilde{w}_\e (t,\gamma(t))\vert\leq C \vert \tilde{y}_\e(t) \vert +  C \vert  D^2 \tilde{w}_\e \vert+  C|\nabla \tilde{w}_\e|+C\e,$$
then $$ \|D^2 \tilde{w}_\e(t,\cdot)\|_{L^\infty(\R^d)} \leq C\int_0^t |y|(s)\mathrm{d}s + C\int_0^t \|D^2 \tilde{w}_\e(s,\cdot)\|_{L^\infty(\R^d)} \mathrm{d}s+  C\int_0^t \|\nabla \tilde{w}_\e(s,\cdot)\|_{L^\infty(\R^d)} \mathrm{d}s +C\e.$$
According to the ODE, we have $$ \frac{\mathrm{d}}{\mathrm{d}t} \vert \tilde{y}_\e\vert \leq C \vert\tilde{y}_\e\vert + \vert C D^2 \tilde{w}_\e (t,\overline{x}(t)) \vert +C\e, $$
and then
$$ |\tilde{y}_\e(t)| \leq C \int_0^t |y(s)| \mathrm{d}s + C \int_0^t \|D^2 \tilde{w}_\e(s,\cdot)\|_{L^\infty(\R^d)}\mathrm{d}s+C\e.$$
We sum the last equations to obtain
\begin{align*}
    &|\tilde{y}_\e(t)| +\|\tilde{w}_\e(t,\cdot)\|_{L^\infty(\R^d)}+ \|\nabla \tilde{w}_\e(t,\cdot)\|_{L^\infty(\R^d)} + \|D^2 \tilde{w}_\e(t,\cdot)\|_{L^\infty(\R^d)} \\ & \leq  \ C \int_0^t\Big( |\tilde{y}_\e(s)| +\|\tilde{w}_\e(s,\cdot)\|_{L^\infty(\R^d)}  + \|\nabla \tilde{w}_\e(s,\cdot)\|_{L^\infty(\R^d)} + \|D^2 \tilde{w}_\e(s,\cdot)\|_{L^\infty(\R^d)} \Big)\mathrm{d}s + C\e.
    \end{align*}
We use the Grönwall lemma to conclude that $(y_\e,w_\e)$ converges in $L_t^\infty W_x^{2,\infty}([0,T]\times \R^d)$ to $(y,w)$.\\\\
\textbf{(ii)} To prove the uniqueness, we use the same method. Let $(y_1,w_1)$ and $(y_2,w_2)$ be two couples of solutions.\\
We define $ y := y_1-y_2$, and $w := w_1-w_2$, they satisfy: 
$$\begin{cases}
    \dot y(t) &= A(t) y(t) + B(t) D^2 w(t,\overline x(t)) C(t)+E(t)y(t)F(t)\\
    \partial_t w &= 2 \nabla u \cdot \nabla w +G(t,x)y(t).
\end{cases}$$
We can determine $w$ along the characteristics that satisfy:
$$\dot{\gamma}(t) = -2 \nabla u (t,\gamma(t)).$$
According to the Cauchy-Lipshitz Theorem, for $t_1 \in [0,T]$, and for $x_1\in \R^d$, there exists a unique trajectory $\gamma$ defined in $[0,t_1]$ such that $\gamma(t_1)=x_1$.
As before, we obtain
$$ \forall t \in [0,t_1], \quad\frac{\mathrm{d}}{\mathrm{d}t}\vert w\vert (t,\gamma(t)) \leq C \vert y\vert(t).$$
We integrate, to find
$$ \forall t\in [0,T], \quad \forall x \in \R^d, \quad |w(t,x)|\leq C\int_0^t |y|(s)\mathrm{d}s. $$
We take the supremum in space of the last inequality, and obtain
$$ \|w(t,\cdot)\|_{L^\infty(\R^d)} \leq C\int_0^t |y|(s)\mathrm{d}s.$$
We differentiate the second equation with respect to $x$, and we use that $\nabla G$ and $D^2 u$ are bounded to obtain the inequality
$$ \frac{\mathrm{d}}{\mathrm{d}t} \vert \nabla w \vert (t,\gamma(t)) \leq C \vert \nabla w \vert (t,\gamma(t)) + C\vert y(t) \vert,$$
then we integrate and  we use that $D^2G$ and $D^3 u$ are bounded to find
$$ \|\nabla w(t,\cdot)\|_{L^\infty(\R^d)} \leq C\int_0^t |y|(s)\mathrm{d}s + C\int_0^t \|\nabla w(s,\cdot)\|_{L^\infty(\R^d)} \mathrm{d}s.$$
Similarly, we obtain an estimate of the second derivative
$$\frac{\mathrm{d}}{\mathrm{d}t} \vert D^2 w (t,\gamma(t))\vert\leq C \vert y(t) \vert +  C \vert  D^2 w \vert+  C|\nabla w| ,$$
then $$ \|D^2 w(t,\cdot)\|_{L^\infty(\R^d)} \leq C\int_0^t |y|(s)\mathrm{d}s + C\int_0^t \|D^2 w(s,\cdot)\|_{L^\infty(\R^d)} \mathrm{d}s+  C\int_0^t \|\nabla w(s,\cdot)\|_{L^\infty(\R^d)} \mathrm{d}s.$$
According to the ODE, we have $$ \frac{\mathrm{d}}{\mathrm{d}t} \vert y\vert \leq C \vert y\vert + \vert C D^2 w (t,\overline{x}(t)) \vert, $$
and then
$$ |y(t)| \leq C \int_0^t |y(s)| \mathrm{d}s + C \int_0^t \|D^2 w(s,\cdot)\|_{L^\infty(\R^d)} \mathrm{d}s.$$
We sum the last equations to obtain
\begin{align*} &|y(t)| +\|w(t,\cdot)\|_{L^\infty(\R^d)}+ \|\nabla w(t,\cdot)\|_{L^\infty(\R^d)} + \|D^2 w(t,\cdot)\|_{L^\infty(\R^d)}\\& \leq C \int_0^t |y(s)| +\|w(s,\cdot)\|_{L^\infty(\R^d)}+ \|\nabla w(s,\cdot)\|_{L^\infty(\R^d)}+ \|D^2 w(s,\cdot)\|_{L^\infty(\R^d)}\mathrm{d}s \end{align*}
We establish the uniqueness using Grönwall's Lemma. Hence the sequences $(y_\e)_\e$ and $(w_\e)_\e$ converges to the unique solution $(y,w)$ of the system.
\begin{flushright}
    $\square$
\end{flushright}
\subsection{\texorpdfstring{Proof of Theorem \ref{thm:approx}}{}}
We proved in the last section that we have for $t\in [0,T],$ $$| I_\e(t) - \mathcal{I}\big(x_\e(t)\big) - K(t)\e |\leq C \e^2.$$
Using the expansion of $x_\e$, we obtain for $t\in [0,T],$
$$|I_\e(t) - I(t) - \e J(t)|\leq  C\e^2.$$
Moreover, we showed that $$\| v_\e - u - \e w\|_{L^\infty_t W^{2,\infty}_x([0,T]\times \R^d)} \leq C\e^2.$$
Hence, $$\| u_\e - u - \e w - \e \log \Big(\frac{r}{\e^{\f d 2}}\Big)\|_{L^\infty_t W^{2,\infty}_x([0,T]\times \R^d)} \leq C\e^2.$$
\begin{flushright}
    $\square$
\end{flushright}
\textbf{Acknowledgments:} The authors would like to thank Sébastien Lion for valuable discussions. This work was supported by the French National Research Agency (ANR) through projects DEEV ANR-20-CE40-0011-01 and ReaCh ANR-23-CE40-0023-01.

\appendix
\section{Proof of Lemma \ref{lemmeK}}
\label{A_K}
We prove this result for $d=1$. Similar arguments can be used to prove the result for $d\geq 2.$ \\ We recall that $J_\e = \log I_\e$ and start from the ODE \fer{ODE-J_e} again
\begin{align*}
\e \frac{\mathrm{d} }{\mathrm{d}t}J_\e(t) = & \ \e^2 \frac{\int_\R n_\e(t,x) \Delta \psi (x)\mathrm{d}x}{\int_\R n_\e(t,x) \psi (x)\mathrm{d}x} + R\big(x_\e(t),I_\e(t)\big) \\& + \frac{\int_\R e^{\frac{u_\e(t,x)-u_\e(t,x_\e(t))}{\e}}\Big(R\big(x,I_\e(t)\big)- R\big(x_\e(t),I_\e(t)\big)\Big) \psi(x)\mathrm{d}x}{\int_\R e^{\frac{u_\e(t,x)-u_\e(t,x_\e(t))}{\e}} \psi (x)\mathrm{d}x}.\nonumber
\end{align*}
According to the estimate \fer{laplace}, we write 
$$\int_\R n_\e(t,x) \Delta \psi (x)\mathrm{d}x = e^{\frac{u_\e(t,x_\e(t))}{\e}} \Delta \psi\big(x_\e(t)\big) \sqrt{\frac{2\pi\e}{|D^2u_\e(t,x_\e(t))|}}+o(\sqrt\e),$$ and 
$$ \int_\R n_\e(t,x) \psi (x)\mathrm{d}x = e^{\frac{u_\e(t,x_\e(t))}{\e}} \psi\big(x_\e(t)\big) \sqrt{\frac{2\pi\e}{|D^2u_\e(t,x_\e(t))|}}+o(\sqrt\e).$$
Thus, we can estimate the first term of the right-hand side in \fer{ODE-J_e} $$\e^2 \frac{\int_\R n_\e(t,x) \Delta \psi (x)\mathrm{d}x}{\int_\R n_\e(t,x) \psi (x)\mathrm{d}x} = \e^2 \frac{\Delta \psi(x_\e(t))}{\psi(x_\e(t))}+ o(\e^2).$$
We use the estimate $ x_\e(t) = \overline{x}(t) + O(\e),$ given in proposition \ref{item1} 
and we obtain: $$\e^2 \frac{\int_\R n_\e(t,x) \Delta \psi (x)\mathrm{d}x}{\int_\R n_\e(t,x) \psi (x)\mathrm{d}x} = \e^2 \frac{\Delta \psi\big(\overline{x}(t)\big)}{\psi\big(\overline{x}(t)\big)}+ o(\e^2).$$
Plugging the latter result in \fer{ODE-J_e}, we find
   \begin{align} 
   \label{Je2}
   \e \frac{\mathrm{d} }{\mathrm{d}t}J_\e(t) =& \  \e^2 \frac{\Delta \psi\big(\overline{x}(t)\big)}{\psi\big(\overline{x}(t)\big)} + R\big(x_\e(t),I_\e(t)\big)\\ &+ \frac{\int_\R e^{\frac{u_\e(t,x)-u_\e(t,x_\e(t))}{\e}}\Big(R\big(x,I_\e(t)\big)- R\big(x_\e(t),I_\e(t)\big)\Big) \psi(x)\mathrm{d}x}{\int_\R  e^{\frac{u_\e(t,x)-u_\e(t,x_\e(t))}{\e}} \psi (x)\mathrm{d}x} +o(\e^2).\nonumber \end{align}
We now look for an expansion of order $O(\e^2)$ for the last term of the right-hand side in \fer{Je2}: $$\frac{\int_\R e^{\frac{u_\e(t,x)-u_\e(t,x_\e(t))}{\e}}\Big(R\big(x,I_\e(t)\big)- R\big(x_\e(t),I_\e(t)\big)\Big) \psi(x)\mathrm{d}x}{\int_\R e^{\frac{u_\e(t,x)-u_\e(t,x_\e(t))}{\e}} \psi (x)\mathrm{d}x}.$$
\\\\ We already have an estimate for the denominator in \fer{laplace} $$ {\int_{\R} e^{\frac{u_\e(t,x)-u_\e(t,x_\e(t))}{\e}} \psi (x)\mathrm{d}x} = \sqrt{\frac{2\pi \e}{\vert D^2u(t,\overline{x}(t))\vert}} \psi \big(\overline{x}(t)\big)+ o(\sqrt{\e}).$$
We now focus on the numerator, we want an expansion of order $O(\e^2\sqrt \e).$ Let $\alpha$ be nonnegative.\\
We split the integral in two parts, and define:
\begin{align*}
    &A_1 :=\int_{|x-x_\e|\leq  \e^{\alpha} } e^{\f {u_\e(x,t)-u_\e(x_\e(t),t)}{\e}} \psi(x) \Big( R\big(x,I_\e(t)\big) - R\big(x_\e(t), I_\e(t)\big) \Big)\mathrm{d}x, \\
    &A_2 := \int_{|x-x_\e|> \e^{\alpha} }  e^{\f {u_\e(x,t)-u_\e(x_\e(t),t)}{\e}} \psi(x) \Big( R\big(x,I_\e(t)\big) - R\big(x_\e(t), I_\e(t)\big) \Big)\mathrm{d}x.
\end{align*}\\
We first show that $A_2 = O(\e^2 \sqrt \e).$ As before, there exists $z,z' \in e(x_\e(t),x)$ such that 
\begin{align*}
     A_2  =&  \int_{\vert x-x_\e(t) \vert > \e^{\alpha}} e^{\frac{D^2 u_\e(t,z)(x-x_\e(t))^2}{2\e}}\Big(\nabla R\big(x_\e(t)),I_\e(t)\big)(x-x_\e(t))  \\&+\frac{1}{2} D^2R\big(z',I_\e(t)\big)(x-x_\e(t))^2 \Big) \psi(x) \mathrm{d}x.
\end{align*}
Using the assumptions \eqref{asr} \fer{asrD2} on $R$, and the concavity estimate \fer{eq.D2u} on $u_\e$ we have:
\begin{align*}
   &\Big \vert \int_{\vert x-x_\e(t) \vert > \e^{\alpha}} e^{\frac{D^2 u_\e(t,z)(x-x_\e(t))^2}{2\e}}\nabla R\big(x_\e(t),I_\e(t)\big)(x-x_\e(t))  \mathrm{d}x\Big\vert \\  & \leq  \vert \nabla R\big(x_\e(t),I_\e(t)\big) \vert \int_{\vert y\vert > \e^{\alpha}} e^{\frac{-\overline{M}_1 }{\e}y^2} \vert y  \vert \mathrm{d}y \\
    &  \leq C \e e^{-\overline{M}_1 \e^{2\alpha-1}} = O(\e^2 \sqrt \e), \text{ for $0<\alpha < \f 12$,}
\end{align*}
and 
\begin{align*}
   &\Big\vert \int_{\vert x-x_\e(t) \vert > \e^{\alpha}} e^{\frac{D^2 u_\e(t,z)(x-x_\e(t))^2}{2\e}}\frac{1}{2} D^2R\big(z',I_\e(t)\big)(x-x_\e(t))^2  \psi(x) \mathrm{d}x \Big\vert\\ 
   & \leq    \|\psi\|_\infty \| D^2R \|_\infty  \e \sqrt \e \int_{\e^{\alpha-\f12}}^{+\infty} e^{-\overline{M}_1y^2}y^2 \mathrm{d}y\\
    &\leq  C \e \sqrt{\e} \int_{\e^{\alpha-\f12}}^{+\infty}e^{-\overline{M}_1 y^2}y^2 \mathrm{d}y = O(\e^2 \sqrt{\e}), \text{ for $0<\alpha < \f 12$.}
\end{align*}
Hence, $A_2 = O(\e^2 \sqrt{\e}).$\\\\
We now expand $A_1$. We use Taylor expansions to write
\begin{align*}
    A_1 =& \int_{\vert x-x_\e(t) \vert \leq \e^\alpha } e^{\frac{\f 12 D^2 u_\e(t,x_\e(t))(x-x_\e(t))^2 + \f 16 D^3 u_\e(t,x_\e(t))(x-x_\e(t))^3 +\frac{1}{24}D^4 u_\e(t,x_\e(t))(x-x_\e(t))^4 + \frac{1}{120} D^5 u_\e(t,z'')(x-x_\e(t))^5}{\e}}\\
   & \times\Big(\psi \big(x_\e(t)\big) + \psi'\big(x_\e(t)\big)\cdot (x-x_\e(t)) +\f 12  \psi''(x_\e(t))(x-x_\e(t))^2 + \f 16 \psi'''(z''')(x-x_\e(t))^3  \Big) \\
   &\times \Big(\nabla R\big(x_\e(t),I_\e(t)\big) (x-x_\e(t))+ \f 12 D^2 R\big(x_\e(t),I_\e(t)\big) (x-x_\e(t))^2 + \f 16 D^3 R\big(z',I_\e(t)\big)(x-x_\e(t))^3 \Big) \mathrm{d}x.
\end{align*}
We make the change of variables $y = \frac{x-x_\e(t)}{\sqrt \e}$ and we obtain
\begin{align*}
    A_1 =&\int_{\vert y \vert \leq \e^{\alpha - \f 12} } e^{\f 12 D^2 u_\e(t,x_\e(t))y^2 + \frac{\sqrt \e}{6} D^3 u_\e(t,x_\e(t))y^3 +\frac{\e}{24}  D^4 u_\e(t,x_\e(t))y^4+\frac{\e\sqrt{\e}}{120} D^5 u_\e(t,z'')(x-x_\e(t))^5} \\
    &\times \Big(\psi (x_\e(t)) + \sqrt \e  \psi'\big(x_\e(t)\big)\cdot y + \f 12 \e \psi''(x_\e(t))y^2 + \f16 \e \sqrt{\e} \psi'''(z''')y^3 \Big) \\
    &\times \Big(\sqrt \e \nabla R\big(x_\e(t),I_\e(t)\big) \cdot y +   \f 12 \e D^2 R\big(x_\e(t),I_\e(t)\big) y^2 + \f 16 \e \sqrt{\e }D^3 R\big(z',I_\e(t)\big)y^3 \Big)\sqrt \e \mathrm{d}y.
\end{align*}
We develop the exponential for $\f 12 >\alpha > \frac{1}{3}$ and we find 
\begin{align*}
    A_1 =& \ \int_{\vert y \vert \leq \e^{\alpha - \f 12} }  e^{\f 12 D^2 u_\e(t,x_\e(t))y^2} \Big(1 + \frac{\sqrt{\e}}{6}D^3 u_\e(t,x_\e(t)) y^3  + \frac{\e}{24} D^4 u_\e(t,x_\e(t))y^4 \\ & + \frac{\e \sqrt{\e}}{120} D^5 u_\e(t,x_\e(t))y^5 + \frac{\e}{72}(D^3 u_\e(t,x_\e(t)))^2 y^6 + O(\e\sqrt{\e}y^7) \Big) \\ &\times \Big(\psi (x_\e(t)) + \sqrt \e  \psi'\big(x_\e(t)\big)\cdot y + \f 12 \e \psi''(x_\e(t))y^2 + \f16 \e \sqrt{\e} \psi'''(z''')y^3 \Big) \\
    &\times \Big(\sqrt \e \nabla R\big(x_\e(t),I_\e(t)\big) \cdot y +   \f 12 \e D^2 R\big(x_\e(t),I_\e(t)\big) y^2 + \f 16 \e \sqrt{\e }D^3 R\big(z',I_\e(t)\big)y^3 \Big)\sqrt \e \mathrm{d}y.
\end{align*}
We develop the latter equation and we keep only the term of order smaller than $O(\e^2).$ We find 
\begin{align*}
    A_1 = &\  \e \psi\big(x_\e(t)\big) \nabla R\big(x_\e(t),I_\e(t)\big)\int_{|y|\leq \e^{\alpha-\f12}} e^{\f12 D^2 u_\e(t,x_\e(t))y^2} y \mathrm{d}y  \\& + \frac{\e \sqrt \e}{2} \psi\big(x_\e(t)\big) D^2 R\big(x_\e(t),I_\e(t)\big) \int_{|y|\leq \e^{\alpha-\f12}} e^{\f12 D^2 u_\e(t,x_\e(t))y^2}  y^2 \mathrm{d}y\\
 &+\frac{\e \sqrt \e}{6} \psi\big(x_\e(t)\big) \nabla R\big(x_\e(t),I_\e(t)\big) D^3 u_\e(t,x_\e(t))\int_{|y|\leq \e^{\alpha-\f12}} e^{\f12 D^2 u_\e(t,x_\e(t))y^2}  y^4 \mathrm{d}y \\
 &+ \e \sqrt{\e} \psi'\big(x_\e(t)\big) \nabla R\big(x_\e(t),I_\e(t)\big)  \int_{|y|\leq \e^{\alpha-\f12}} e^{\f12 D^2 u_\e(t,x_\e(t))y^2} y^2 \mathrm{d}y\\&
 + \frac{\e^2}{24} \psi\big(x_\e(t)\big) \nabla R\big(x_\e(t),I_\e(t)\big)  D^4 u_\e(t,x_\e(t)) \int_{|y|\leq \e^{\alpha-\f12}}  e^{\f12 D^2 u_\e(t,x_\e(t))y^2}  y^5 \mathrm{d}y \\
 & + \frac{\e^2}{72} \psi\big(x_\e(t)\big) \nabla R\big(x_\e(t),I_\e(t)\big) (D^3 u_\e(t,x_\e(t)))^2 \int_{|y|\leq \e^{\alpha-\f12}}  e^{\f12 D^2 u_\e(t,x_\e(t))y^2}  y^7 \mathrm{d}y \\
 & + \frac{\e^2}{6} \psi'\big(x_\e(t)\big) \nabla R\big(x_\e(t),I_\e(t)\big) D^3 u_\e(t,x_\e(t)) \int_{|y|\leq \e^{\alpha-\f12}}  e^{\f12 D^2 u_\e(t,x_\e(t))y^2}  y^5 \mathrm{d}y\\
 &+ \frac{\e^2}{12} \psi(x_\e(t)) D^2 R\big(x_\e(t),I_\e(t)\big) D^3 u_\e(t,x_\e(t)) \int_{|y|\leq \e^{\alpha-\f12}}  e^{\f12 D^2 u_\e(t,x_\e(t))y^2}  y^5 \mathrm{d}y\\
 & + \frac{\e^2}{6} \psi(x_\e(t)) D^3 R(x_\e(t),I_\e(t)) \int_{|y|\leq \e^{\alpha-\f12}}  e^{\f12 D^2 u_\e(t,x_\e(t))y^2}  y^3 \mathrm{d}y\\
 & + \frac{\e^2}{2} \psi'(x_\e(t)) D^2 R(x_\e(t),I_\e(t) \int_{|y|\leq \e^{\alpha-\f12}}  e^{\f12 D^2 u_\e(t,x_\e(t))y^2}  y^3 \mathrm{d}y \\
 & + \frac{\e^2}{2} \psi''(x_\e(t)) \nabla R\big(x_\e(t),I_\e(t)\big) \int_{|y|\leq \e^{\alpha-\f12}}  e^{\f12 D^2 u_\e(t,x_\e(t))y^2}  y^3 \mathrm{d}y + O(\e^2 \sqrt{\e}).
\end{align*}
We notice that every term of order $\e^2$ is equal to $0$, because we integrate an odd function over a symmetric interval. We now need to compute the remaining integrals. To this end, we use the result of Lemma \ref{lemma2} to find
\begin{align*}
    \int_{|y|\leq \e^{\alpha-\f12}} e^{\f12 D^2 u_\e(t,x_\e(t))y^2}  y^2 \mathrm{d}y = & \ \sqrt{2 \pi} |D^2 u_\e (t,\bar x (t))|^{-\frac{3}{2}} +O(\e) =  \sqrt{2 \pi} |D^2 u (t,\bar x (t))|^{-\frac{3}{2}} +O(\e), \\
    \int_{|y|\leq \e^{\alpha-\f12}} e^{\f12 D^2 u_\e(t,x_\e(t))y^2}  y^4 \mathrm{d}y = & \ 3\sqrt{2 \pi} |D^2 u_\e(t,\bar x (t))|^{-\frac{5}{2}} +O(\e)= 3\sqrt{2 \pi} |D^2 u(t,\bar x (t))|^{-\frac{5}{2}} +O(\e).
\end{align*}
Finally, we proved that
 $A_1 = \e\sqrt{\e} \bar{f}(t) + O(\e^2\sqrt{\e})$ for $t$ in $[0,T]$ with 
\begin{align*}
    \bar{f}(t) := & \ \sqrt{2\pi} \nabla \psi \big(\overline{x}(t)\big) \cdot \nabla R\big(\overline{x}(t),I(t)\big)  \vert D^2 u(t,\overline{x}(t)) \vert^{-\frac{3}{2}} + \frac{\sqrt{2\pi}}{2 } \psi \big(\overline{x}(t)\big) D^2 R\big(\overline{x}(t),I(t)\big)\vert D^2 u(t,\overline{x}(t)) \vert^{-\frac{3}{2}}\\
    & + \frac{ \sqrt{2\pi}}{2} \psi\big(\overline{x}(t)\big) \nabla R\big(\overline{x}(t),I(t)\big)D^3 u(t,\overline{x}(t)) \vert D^2 u(t,\overline{x}(t)) \vert^{-\frac{5}{2}}.
\end{align*}
We combine the latter result with the equivalent of the denominator, and we obtain: $$ \e \f{\mathrm{d} }{\mathrm{d} t }J_\e(t)= R\big(x_\e(t),I_\e(t)\big)+\e f(t) +O(\e^2), $$
with 
\begin{align*}
    f(t) =& \Big[\frac{\nabla \psi \big(\overline{x}(t)\big)}{ \psi\big(\overline{x}(t)\big)} 
    \cdot \nabla R\big(\overline{x}(t),I(t)\big)   + \frac{1}{2 } D^2 R\big(\overline{x}(t),I(t)\big)\Big]\vert D^2 u(t,\overline{x}(t)) \vert^{-1}\\ 
    &+  \frac{1}{2}  \nabla R\big(\overline{x}(t),I(t)\big)D^3 u(t,\overline{x}(t)) \vert D^2 u(t,\overline{x}(t)) \vert^{-2}.
\end{align*}
Moreover, we have $R \left( x_\e(t), \mathcal I(x_\e(t) \right) =0.$
We differentiate with respect to time, and we use the error estimate \fer{errorx} on $x_\e$,  from Proposition \ref{item1}, to obtain 
$$
\f{\mathrm{d}}{\mathrm{d}t} \mathcal J\big(x_\e(t)\big) = - \f{\nabla_x R\big(\overline x(t),I(t)\big) \left(-D^2 u(t,\overline x(t)) \right)^{-1}\nabla R\big(\overline x(t),I(t)\big)}{I(t) \f{\p }{\p I}R\big(\overline x(t),I(t)\big)}+O(\e)=g(t)+O(\e).
$$
We combine the two ODE to write:
$$
 \e \f{\mathrm{d} }{\mathrm{d} t }\big(J_\e(t)-\mathcal J\big(x_\e(t)\big)\big)= R\big(x_\e(t),I_\e(t)\big) - R\big(x_\e(t), \mathcal I\big(x_\e(t)\big)\big)+\e h(t) +O(\e^2), \quad \text{with $h(t)=f(t)-g(t)$}.$$
We define $ k_\e(t) : = \frac{J_\e(t)-\mathcal J(x_\e(t))}{\e},$ and obtain the following ODE
$$\e \f{\mathrm{d} }{\mathrm{d} t }k_\e(t) = \f{R\big(x_\e(t),I_\e(t)\big) - R\big(x_\e(t), \mathcal I\big(x_\e(t)\big)\big)}{\e}+ h(t) +O(\e).$$
Using the mean value theorem with $\tilde{R}(\cdot,J) = R(\cdot, e^J),$ we obtain
$$ \e \f{\mathrm{d} }{\mathrm{d} t }k_\e(t) =\widetilde I_\e(t)\, \f{\p}{\p I} R\big(x_\e(t),\widetilde I_\e(t)\big) \, k_\e(t) +h(t)+O(\e), \quad \text{with $\widetilde I_\e(t)$ between $I_\e(t)$ and $\mathcal I(x_\e(t))$}.$$
Let us show that $(k_\e)_\e$ converges when $\e$ approaches $0$. The application $k_\e$ satisfies an ODE of order $1$, which can be solved
$$k_\e(t) = k_\e(0) e^{\frac{1}{\e}\int_{0}^t \tilde I _\e(s) \frac{\partial}{\p I} R(x_\e(s)),\tilde I _\e(s)) \mathrm{d}s} + k_{p,\e}(t),$$ where $k_{p,\e}$ is a particular solution.
We determine $k_{p,\e}$ using the method of variation of constants. We write $ k_{p,\e}(t) := \lambda_\e (t) e^{\frac{1}{\e}\int_{0}^t \tilde I _\e(t) \frac{\partial}{\p I} R(x_\e(s)),\tilde I _\e(t))\mathrm{d}s}.$ We differentiate to obtain: 
$$\e \f{\mathrm{d} }{\mathrm{d} t }k_{p,\e}(t) = \big(\e \lambda'_\e(t) + \lambda_\e\tilde I _\e(t) \frac{\partial}{\p I} R\big(x_\e(t),\tilde I _\e(t)\big) \big) \exp \Big(\frac{1}{\e}\int_{0}^t \tilde I _\e(s) \mathrm{d}s\frac{\partial}{\p I} R\big(x_\e(s),\tilde I _\e(s)\big) \mathrm{d}s\Big),$$
and we use the ODE to find: 
$$ k_{p,\e}(t) = \frac{1}{\e} \int_0^t  e^{\frac{1}{\e}\int_{u}^t \tilde I _\e(s) \frac{\partial}{\p I} R(x_\e(s),\tilde I _\e(s))\mathrm{d}s}(h(u)+O(\e))\mathrm{d}u 
.$$ We write the latter integral as follows
$$ k_{p,\e}(t) = A_1(t)+A_2(t),$$ with 
$$ A_1(t) =  \frac{1}{\e} \int_0^t  e^{\frac{1}{\e}\int_{u}^t \tilde I _\e(s) \frac{\partial}{\p I} R(x_\e(s),\tilde I _\e(s))\mathrm{d}s}h(u) \mathrm{d}u, \ A_2(t) = \int_0^t e^{\frac{1}{\e}\int_{u}^t \tilde I _\e(s) \frac{\partial}{\p I} R(x_\e(s),\tilde I _\e(s))\mathrm{d}s}O(1) \mathrm{d}u.$$
We first estimate $A_2.$ Thanks to \fer{asrDi}, we have that $\f{\p}{\p I} R\leq -K_2.$ Moreover, $\tilde I_\e(t) \in e(I_\e(t),\mathcal{I}(x_\e(t))),$ thanks to Theorem \ref{t2.1} and Proposition \ref{item1}, we deduce that $\tilde I_\e(t) \geq \frac{I_m}{2}>0.$ We obtain
$$ |A_2(t)| \leq C\int_0^t e^{-\frac{(t-u)}{\e}\cdot\frac{I_m \overline K_2}{2}} \mathrm{d}u = \frac{2C\e}{I_m \overline K_2}(1-e^{-\frac{t}{\e}\cdot \frac{I_m \overline K_2}{2}}) = O(\e).$$
We will now compute $A_1$ using integration by parts. We find
\beq \label{A1} A_1(t) = \Big[ -\frac{h(u)}{\tilde I_\e(u) \partial_I R(x_\e(u),I_\e(u))}e^{\frac{1}{\e}\int_{u}^t \tilde I _\e(s) \frac{\partial}{\p I} R(x_\e(s),\tilde I _\e(s))\mathrm{d}s}\Big]_0^t + \int_0^t h'(u) e^{\frac{1}{\e}\int_{u}^t \tilde I _\e(s) \frac{\partial}{\p I} R(x_\e(s),\tilde I _\e(s))\mathrm{d}s} \mathrm{d}u.\eeq
Thanks to \fer{asrDi}, Theorem \ref{t2.1} and Proposition \ref{item1}, the integral of the right-hand side of \fer{A1} is of order $O(\e).$ We deduce that
$$A_1(t) = K(t) - K(0)e^{\frac{1}{\e}\int_{0}^t \tilde I _\e(s) \frac{\partial}{\p I} R(x_\e(s),\tilde I _\e(s))\mathrm{d}s}+O(\e),$$
with $$ K(t) = -\frac{h(t)}{I(t)\frac{\partial}{\partial I} R(\overline{x}(t),I(t))}.$$ Finally, for $t>0$, we have $k_{p,\e}(t) = K(t) +O(\e)$ and $k_{p,\e}(0) = O(\e).$ We conclude that $k_\e(t) = K(t) +O(\e)$ for all $t>0$. 
Consequently, we deduce that $ J_\e(t) = \mathcal{J}\big(x_\e(t)\big) + \e K(t) + O(\e^2)$, and $ I_\e(t) = \mathcal{I}\big(x_\e(t)\big)(1+ \e K(t) +O(\e^2))$. This concludes the proof of Lemma \ref{lemmeK}. 
\section{Approximation of the moments -- Proof of Theorem \ref{thm:moments}}
\label{A_mom}
We prove the result for $d=1.$ The proof could be adapted for $d\geq 2.$ 
\\\\ \textbf{(i) Proof of \fer{m1}: the expansion of $M_{1,\e}$.}\\
    We use the asymptotic expansion from Theorem \ref{thm:approx} to write
    \begin{align*}
    M_{1,\e}(t) &= \frac{1}{I_\e(t)} \int x e^{\frac{u_\e(t,x)}{\e}}\mathrm{d}x\\
    &= \overline{x}(t) + \frac{1}{I_\e(t)} \int (x -\overline{x}(t))e^{\frac{u_\e(t,x)}{\e}}\mathrm{d}x\\
   & = \overline{x}(t) + \frac{\int (x -\overline{x}(t))e^{\frac{u(t,x)-u(t,\overline{x}(t))}{\e}+ v(t,x) + O(\e)}\mathrm{d}x}{\int e^{\frac{u(t,x)-u(t,\overline{x}(t))}{\e}+ v(t,x) + O(\e)} \mathrm{d}x}.
    \end{align*}
We proceed as in the proof of Lemma \ref{lemme1}, equation \fer{laplace}, to find an equivalent of the denominator of the last term of the right-hand side $$ {\int e^{\frac{u(t,x)-u(t,\overline{x}(t))}{\e}+ v(t,x) + O(\e)}dx} = \sqrt{\frac{2\pi\e}{|D^2 u(t,\overline{x}(t))|}}\exp(v(t,\overline{x}(t))) + o(\sqrt{\e}).$$
Now, we need to compute $\int (x-\overline{x}(t))e^{\frac{u(t,x)-u(t,\overline{x}(t))}{\e}+ v(t,x)+O(\e)}\mathrm{d}x.$\\
Let $\alpha$ be a nonnegative constant.
We split the integral in two terms: 
\begin{align*}
    A_1 &:= \int_{ |x-\overline{x}(t)| \leq \e^\alpha} (x -\overline{x}(t))e^{\frac{u(t,x)-u(t,\overline{x}(t))}{\e}+ v(t,x)+O(\e)}\mathrm{d}x,\\
    A_2 &:= \int_{ |x-\overline{x}(t)| > \e^\alpha} (x -\overline{x}(t))e^{\frac{u(t,x)-u(t,\overline{x}(t))}{\e}+ v(t,x) +O(\e)}\mathrm{d}x.
\end{align*}
We know that $\overline{x}(t)$ maximizes $u(t,\cdot)$, we write the following Taylor expansion
$$u(t,x) = u(t,\overline{x}(t)) + \f 12 D^2 u(t,\overline{x}(t))(x-\overline{x}(t))^2 + \f 16 D^3 u(t,\overline{x}(t))(x-\overline{x}(t))^3+ \f{1}{24} D^4 u(t,x')(x-\overline{x}(t))^4,$$
$\text{with } x' \in e(x,\overline{x}(t)).$
Thus, we deduce that
\begin{align*}
    A_1 =  \int_{ |x-\overline{x}(t)| \leq \e^\alpha} (x -\overline{x}(t))e^{\f{1}{2\e} D^2u(t,\overline{x}(t))(x-\overline{x}(t))^2 + \f 16 D^3 u(t,\overline{x}(t))(x-\overline{x}(t))^3+ \f{1}{24} D^4 u(t,x')(x-\overline{x}(t))^4 v(t,x)+O(\e)}
    \mathrm{d}x.
\end{align*}
We make the following change of variable $y = \frac{x-\overline{x}(t)}{\sqrt \e}$ and for $\f 1 3 < \alpha < \f 1 2,$ and we obtain
\begin{align*}
    A_1 =& \ \e \int_{ |y| \leq \e^{\alpha-\f 12}} y e^{\f{1}{2} D^2u(t,\overline{x}(t))y^2 + \f{\sqrt \e}{6} D^3 u(t,\overline{x}(t))y^3 +  \f{\e}{24} D^4 u(t,x')y^4+v(t,\overline{x}(t) + \sqrt \e y)+O(\e)}\mathrm{d}y\\
   =& \ \e \int_{ |y| \leq \e^{\alpha-\f 12}} y e^{\f{1}{2} D^2u(t,\overline{x}(t))y^2} \Big(1+\f{\sqrt \e}{6} D^3 u(t,\overline{x}(t))y^3+ \sqrt \e \nabla v(t,\overline{x}(t) )\cdot y + \\ & + O(\e y^6+ \e +\e  y^4+\e y^2) \Big) e^{v(t,\overline{x}(t))} \mathrm{d}y \\
  = & \ \f{\e \sqrt \e}{6} D^3 u(t,\overline{x}(t)) \int_\R e^{\f{1}{2} D^2u(t,\overline{x}(t))y^2}y^4 \mathrm{d} y + \e \sqrt \e \nabla v (t,\overline{x}(t)) \int_\R e^{\f{1}{2} D^2u(t,\overline{x}(t))y^2}y^2 \mathrm{d} y +O(\e^2) \\=& \  \f{\sqrt{2\pi}}{2} \cdot \e \sqrt \e|D^2 u(t,\overline{x}(t))|^{-\frac{5}{2}} D^3 u(t,\overline{x}(t)) + \sqrt{2 \pi} \cdot \e \sqrt \e | D^2 u(t,\overline{x}(t))|^{-\frac{3}{2}} \nabla v (t,\overline{x}(t)) + O(\e^2).
\end{align*}
Indeed we use that $D^4 u$ is bounded and that we have $- 2\underline M_1 \leq D^2 u \leq -2\overline M_1$ thanks to the Theorem \ref{t2.2}.
Now, we prove that $A_2 = O(\e^2).$\\
We compute, for $\e$ small enough 
\begin{align*}
    |A_2| \leq&  \int_{ |x-\overline{x}(t)| > \e^\alpha} |x -\overline{x}(t)|e^{\frac{u(t,x)-u(t,\overline{x}(t))}{\e}+ v(t,x)+O(\e)} \mathrm{d}x\\
     \leq& \ \e \int_{ |y| > \e^{\alpha-\f 12}} |y|e^{-\overline M_1 y^2+ v(t,\overline{x}(t) + \sqrt \e y)+O(\e)}\mathrm{d}y \leq C\e \int_{ |y| > \e^{\alpha-\f 12}} |y|e^{-\overline M_1 y^2}\mathrm{d}y = O(\e^2).
\end{align*}
Therefore, we conclude that for $t\in[0,T]$, $$M_{1,\e}(t) = \overline{x}(t) + \f{\e }{2}  |D^2 u(t,\overline{x}(t))|^{-2} D^3 u(t,\overline{x}(t)) + \e | D^2 u(t,\overline{x}(t))|^{-1}\nabla v (t,\overline{x}(t))  + O(\e \sqrt{\e}).$$
\textbf{(ii) Proof of \fer{m2}: the expansion of $M_{2,\e}^c$. }
\\Using the asymptotic expansion \fer{m1} of $M_{1,\e}$, we can write that
\begin{align*}
    M_{2,\e}^c(t) =& \frac{\int \big(x-\overline{x}(t) +M_1(t)\e +O(\e\sqrt{\e})\big)^2 e^{\frac{u(t,x)-u(t,\overline{x}(t))}{\e}+v(t,x)+O(\e)} \mathrm{d}x}{\sqrt{\frac{2\pi\e}{|D^2u(t,\overline{x}(t))|}}\exp(v(t,\overline{x}(t))+O(\e)}.
\end{align*}
We need an expansion of the numerator.  For $\f 13 < \alpha < \f 12$, we will split the integral into two parts, 
\begin{align*}
    A_1 :=& \int_{ |x-\overline{x}(t)| \leq \e^\alpha} \big(x-\overline{x}(t) +M_1(t)\e +O(\e\sqrt{\e})\big)^2 e^{\frac{u(t,x)-u(t,\overline{x}(t))}{\e}+v(t,x)+O(\e)} \mathrm{d}x,\\
    A_2 :=& \int_{ |x-\overline{x}(t)| > \e^\alpha} \big(x-\overline{x}(t) +M_1(t)\e +O(\e\sqrt{\e})\big)^2 e^{\frac{u(t,x)-u(t,\overline{x}(t))}{\e}+v(t,x)+O(\e)} \mathrm{d}x
    .
\end{align*}
Next, we obtain
\begin{align*}
    A_1 =& \, \sqrt \e \int_{ |y| \leq \e^{\alpha-\f 12}} (\sqrt \e y +M_1(t)\e +O(\e\sqrt{\e}))^2 e^{\f 12 D^2 u(t,\overline{x}(t))y^2}\\ &\times\Big(1+\f{\sqrt \e}{6} D^3 u(t,\overline{z})y^3+ O(\e y^6)+O(\e) + \sqrt \e \nabla v(t,\overline{x}(t) )+O(\e y^2) \Big)\mathrm{d}y \,e^{v(t,\overline{x}(t))} \\
     =& \ \e \sqrt \e \int_\R y^2 e^{\f 12 D^2 u(t,\overline{x}(t))y^2}  \mathrm{d}y \times e^{v(t,\overline{x}(t))}+ O(\e^2 \sqrt \e)\\
    =&\ \sqrt{2\pi}\cdot \e \sqrt \e|D^2u(t,\overline{x}(t))|^{-\f 32} e^{v(t,\overline{x}(t))} + O(\e^2 \sqrt \e).
\end{align*}
With methods similar to those used before, we show that $A_2 = O(\e^2 \sqrt \e)$. Thus, we deduce that $M_{2,\e}^c(t) = \e M_2(t) +O(\e^2)$ with $M_2(t) = |D^2 u(t,\overline{x}(t))|^{-1}$.\\\\
\textbf{(iii) Proof of \fer{m2k}: the expansion of $M_{2k,\e}^c$. }
\\Using the asymptotic expansion \fer{m1} of $M_{1,\e}$, we can write that
\begin{align*}
    M_{2k,\e}^c(t) =& \frac{\int \big(x-\overline{x}(t) +M_1(t)\e +O(\e\sqrt{\e})\big)^{2k} e^{\frac{u(t,x)-u(t,\overline{x}(t))}{\e}+v(t,x)+O(\e)} \mathrm{d}x}{\sqrt{\frac{2\pi\e}{|D^2u(t,\overline{x}(t))|}}\exp(v(t,\overline{x}(t))+O(\e)}.
\end{align*}
We need an expansion of the numerator.  For $\f 13 < \alpha < \f 12$, we will split the integral into two parts, 
\begin{align*}
    A_1 :=& \int_{ |x-\overline{x}(t)| \leq \e^\alpha} \big(x-\overline{x}(t) -M_1(t)\e +O(\e\sqrt{\e})\big)^{2k}  e^{\frac{u(t,x)-u(t,\overline{x}(t))}{\e}+v(t,x)+O(\e)} \mathrm{d}x,\\
    A_2 :=& \int_{ |x-\overline{x}(t)| > \e^\alpha} \big(x-\overline{x}(t) -M_1(t)\e +O(\e\sqrt{\e})\big)^{2k}  e^{\frac{u(t,x)-u(t,\overline{x}(t))}{\e}+v(t,x)+O(\e)} \mathrm{d}x
    .
\end{align*}
Next, we obtain
\begin{align*}
    A_1 =& \, \sqrt \e \int_{ |y| \leq \e^{\alpha-\f 12}} (\sqrt \e y -M_1(t)\e +O(\e\sqrt{\e}))^{2k} e^{\f 12 D^2 u(t,\overline{x}(t))y^2}\\ &\times\Big(1+\f{\sqrt \e}{6} D^3 u(t,\overline{z})y^3+ O(\e y^6)+O(\e) + \sqrt \e \nabla v(t,\overline{x}(t))y+O(\e y^2) \Big)\mathrm{d}y \,e^{v(t,\overline{x}(t))} \\
     =& \ \e^k \sqrt \e \int_\R y^{2k} e^{\f 12 D^2 u(t,\overline{x}(t))y^2}  \mathrm{d}y \times e^{v(t,\overline{x}(t))}+ O(\e^{k+1}).
\end{align*}
We denote $\Gamma_k(t):= \int y^{2k} e^{\frac{1}{2}D^2u(t,\overline{x}(t))y^2}\mathrm{d}y.$ By integration by parts, we obtain
\beq \label{gamma} \Gamma_k(t) = \sqrt{2\pi} \cdot \frac{(2k)!}{2^k k!}M_2(t)^{k+\f 1 2}.\eeq
On the other hand, we find that $A_2 = O(\e^{k+1})$. Thus, we deduce that $M_{2k,\e}^c(t) = \e^{k} M_{2k}(t) +O(\e^k \sqrt{\e})$.\\\\
\textbf{(iv) Proof of \fer{m2k1}: the expansion of $M_{2k+1,\e}^c$. }
\begin{align*}
    M_{2k+1,\e}^c(t) =& \frac{\int \big(x-\overline{x}(t) +M_1(t)\e +O(\e\sqrt{\e})\big)^{2k+1} e^{\frac{u(t,x)-u(t,\overline{x}(t))}{\e}+v(t,x)+O(\e)} \mathrm{d}x}{\sqrt{\frac{2\pi\e}{|D^2u(t,\overline{x}(t))|}}\exp(v(t,\overline{x}(t))+O(\e)}.
\end{align*}
We need an expansion of the numerator.  For $\f 13 < \alpha < \f 12$, we will split the integral into two parts, 
\begin{align*}
    A_1 :=& \int_{ |x-\overline{x}(t)| \leq \e^\alpha} \big(x-\overline{x}(t) -M_1(t)\e +O(\e\sqrt{\e})\big)^{2k+1}  e^{\frac{u(t,x)-u(t,\overline{x}(t))}{\e}+v(t,x)+O(\e)} \mathrm{d}x,\\
    A_2 :=& \int_{ |x-\overline{x}(t)| > \e^\alpha} \big(x-\overline{x}(t) -M_1(t)\e +O(\e\sqrt{\e})\big)^{2k+1}  e^{\frac{u(t,x)-u(t,\overline{x}(t))}{\e}+v(t,x)+O(\e)} \mathrm{d}x
    .
\end{align*}
Next, we simplify $A_1$ and use the expression     \eqref{gamma} of $\Gamma_k(t)$ to obtain
\begin{align*}
    A_1 =& \, \sqrt \e \int_{ |y| \leq \e^{\alpha-\f 12}} (\sqrt \e y -M_1(t)\e +O(\e\sqrt{\e}))^{2k+1} e^{\f 12 D^2 u(t,\overline{x}(t))y^2}\\ &\times\Big(1+\f{\sqrt \e}{6} D^3 u(t,\overline{z})y^3+ O(\e y^6)+O(\e) + \sqrt \e \nabla v(t,\overline{x}(t) )y+O(\e y^2) \Big)\mathrm{d}y \,e^{v(t,\overline{x}(t))} \\
     =& \ 0+ \e^{k+1} \sqrt \e \frac{D^3u(t,\overline{x}(t))}{6}\int_\R y^{2k+4} e^{\f 12 D^2 u(t,\overline{x}(t))y^2}  \mathrm{d}y e^{v(t,\overline{x}(t))} \\ & \ +\e^{k+1} \sqrt \e \nabla v(t,\overline{x}(t))\int_\R y^{2k+2} e^{\f 12 D^2 u(t,\overline{x}(t))y^2}  \mathrm{d}y  e^{v(t,\overline{x}(t))} \\
     & \ - \e^{k+1} \sqrt \e M_1(t)(2k+1)\int_\R y^{2k} e^{\f 12 D^2 u(t,\overline{x}(t))y^2}  \mathrm{d}y  e^{v(t,\overline{x}(t))}+ O(\e^{k+2})\\
     =& \ \e^{k+1} \sqrt \e e^{v(t,\overline{x}(t))} \Big[ \frac{D^3u(t,\overline{x}(t))}{6} \Gamma_{k+2}(t) + \nabla v (t,\overline{x}(t)) \Gamma_{k+1}(t) - (2k+1)M_1(t) \Gamma_k(t)\Big]+ O(\e^{k+2}).
\end{align*}
Separately, using the expression \eqref{M11} of $M_1(t)$, we compute 
\begin{align*}
    (2k+1)M_1(t) \Gamma_k(t) = \Gamma_{k+1}(t)\Big( \frac{1}{2}M_2(t)D^3 u(t,\overline{x}(t))+ \nabla v(t,\overline{x}(t))\Big).
\end{align*}
Finally, we find that
\begin{align*}
    A_1 = & \ \e^{k+1} \sqrt \e e^{v(t,\overline{x}(t))} D^3 u(t,\overline{x}(t)) \Big(\frac{\Gamma_{k+2}(t)}{6}-\frac{1}{2}M_2(t)\Gamma_{k+1}(t)\Big)+ O(\e^{k+2})\\
    = & \ \e^{k+1} \sqrt \e e^{v(t,\overline{x}(t))} D^3 u(t,\overline{x}(t)) \Big( \frac{2k+3}{6} M_2(t) - \f 1 2 M_2(t)\Big) \Gamma_{k+1}(t)+ O(\e^{k+2})\\
    = & \ \e^{k+1} \sqrt \e e^{v(t,\overline{x}(t))} D^3 u(t,\overline{x}(t)) \cdot \frac{k}{3}\cdot \frac{(2(k+1))!}{2^{k+1}(k+1)!}M_2(t)^{k+2 + \f 1 2} \sqrt{2 \pi}+ O(\e^{k+2}). 
\end{align*}
On the other hand, we find that $A_2 = O(\e^{k+2})$. Hence, $M_{2k+1,\e}^c(t) = \e^{k+1} M_{2k+1}(t) +O(\e^{k+1} \sqrt{\e})$. This ends the proof.
\printbibliography
\end{document}